%Paper: q-alg/9512025
%From: Javier Mas <JAMAS@GAES.USC.ES>
%Date: Fri, 22 Dec 1995 22:53:46 WST
%Date (revised): Sat, 23 Dec 1995 05:06:56 WST

%
%%%%%%%%%%%%%%%%%%%%%%%%%%%%%%%%%%%%%%%%%%%%%%%%%%%%%%%%%%%%%%%%
%
%     Title:   The algebra of $q$-pseudodifferential symbols and the
%              $q$-$\WKP^{(n)}$ algebra
%     Authors:  J. Mas and M. Seco
%
%
\nonstopmode
\catcode`\@=11 % this allows for tricky names
%
%%%%%%%%%%%%%%%%%%%%%%%%%%%%%%%%%%%%%%%%%%%%%%%%%%%%%%%%%%%%%%%%%
%
%  First some font definitions
%
\font\seventeenrm=cmr17

\font\twelverm=cmr12
\font\ninerm=cmr9
\font\sixrm=cmr6

\font\seventeenbf=cmbx12 at 17pt
\font\fourteenbf=cmbx12 at 14pt
\font\twelvebf=cmbx12
\font\ninebf=cmbx9
\font\sixbf=cmbx6

\font\seventeeni=cmmi12 at 17pt             \skewchar\seventeeni='177
\font\fourteeni=cmmi12 at 14pt              \skewchar\fourteeni='177
\font\twelvei=cmmi12                        \skewchar\twelvei='177
\font\ninei=cmmi9                           \skewchar\ninei='177
\font\sixi=cmmi6                            \skewchar\sixi='177

\font\seventeensy=cmsy10 scaled\magstep3    \skewchar\seventeensy='60
\font\fourteensy=cmsy10 scaled\magstep2     \skewchar\fourteensy='60
\font\twelvesy=cmsy10 at 12pt               \skewchar\twelvesy='60
\font\ninesy=cmsy9                          \skewchar\ninesy='60
\font\sixsy=cmsy6                           \skewchar\sixsy='60

\font\seventeenex=cmex10 scaled\magstep3
\font\fourteenex=cmex10 scaled\magstep2
\font\twelveex=cmex10 at 12pt

%\font\ninex=cmex9
\font\ninex=cmex10 at 9pt
%\font\sevenex=cmex7
\font\sevenex=cmex10 at 9pt
%\font\sixex=cmex7 at 6pt
\font\sixex=cmex10 at 6pt
%\font\fivex=cmex7 at 5pt
\font\fivex=cmex10 at 5pt

\font\seventeensl=cmsl10 scaled\magstep3
\font\fourteensl=cmsl10 scaled\magstep2
\font\twelvesl=cmsl10 scaled\magstep1
\font\ninesl=cmsl10 at 9pt
\font\sevensl=cmsl10 at 7pt
\font\sixsl=cmsl10 at 6pt
\font\fivesl=cmsl10 at 5pt

\font\seventeenit=cmti12 scaled\magstep2
\font\fourteenit=cmti12 scaled\magstep1
\font\twelveit=cmti12

\font\seventeentt=cmtt12 scaled\magstep2
\font\fourteentt=cmtt12 scaled\magstep1
\font\twelvett=cmtt12

\font\seventeencp=cmcsc10 scaled\magstep3
\font\fourteencp=cmcsc10 scaled\magstep2
\font\twelvecp=cmcsc10 scaled\magstep1
\font\tencp=cmcsc10
%\font\eightcp=cmcsc8

\newfam\cpfam

\font\seventeenss=cmss17
\font\fourteenss=cmss12 at 14pt
\font\twelvess=cmss12
\font\tenss=cmss10
\font\niness=cmss9

\font\sevenss=cmss8 at 7pt
\font\sixss=cmss8 at 6pt
\font\fivess=cmss8 at 5pt
\newfam\ssfam
\newdimen\b@gheight             \b@gheight=12pt
\newcount\f@ntkey               \f@ntkey=0
\def\f@m{\afterassignment\samef@nt\f@ntkey=}
\def\samef@nt{\fam=\f@ntkey \the\textfont\f@ntkey\relax}
\def\rm{\f@m0 }
\def\mit{\f@m1 }         
\def\cal{\f@m2 }
\def\it{\f@m\itfam}
\def\sl{\f@m\slfam}
\def\bf{\f@m\bffam}
\def\tt{\f@m\ttfam}
\def\ssf{\f@m\ssfam}
\def\caps{\f@m\cpfam}
\def\seventeenpoint{\relax
    \textfont0=\seventeenrm          \scriptfont0=\twelverm
      \scriptscriptfont0=\ninerm
    \textfont1=\seventeeni           \scriptfont1=\twelvei
      \scriptscriptfont1=\ninei
    \textfont2=\seventeensy          \scriptfont2=\twelvesy
      \scriptscriptfont2=\ninesy
    \textfont3=\seventeenex          \scriptfont3=\twelveex
      \scriptscriptfont3=\ninex
    \textfont\itfam=\seventeenit    %\scriptfont\itfam=\twelveit
    \textfont\slfam=\seventeensl    %\scriptfont\slfam=\twelvesl
      \scriptscriptfont\slfam=\ninesl
    \textfont\bffam=\seventeenbf     \scriptfont\bffam=\twelvebf
      \scriptscriptfont\bffam=\ninebf
    \textfont\ttfam=\seventeentt
    \textfont\cpfam=\seventeencp
    \textfont\ssfam=\seventeenss     \scriptfont\ssfam=\twelvess
      \scriptscriptfont\ssfam=\niness
    \samef@nt
    \b@gheight=17pt
    \setbox\strutbox=\hbox{\vrule height 0.85\b@gheight
                                depth 0.35\b@gheight width\z@ }}
\def\fourteenpoint{\relax
    \textfont0=\fourteencp          \scriptfont0=\tenrm
      \scriptscriptfont0=\sevenrm
    \textfont1=\fourteeni           \scriptfont1=\teni
      \scriptscriptfont1=\seveni
    \textfont2=\fourteensy          \scriptfont2=\tensy
      \scriptscriptfont2=\sevensy
    \textfont3=\fourteenex          \scriptfont3=\twelveex
      \scriptscriptfont3=\tenex
    \textfont\itfam=\fourteenit     \scriptfont\itfam=\tenit
    \textfont\slfam=\fourteensl     \scriptfont\slfam=\tensl
      \scriptscriptfont\slfam=\sevensl
    \textfont\bffam=\fourteenbf     \scriptfont\bffam=\tenbf
      \scriptscriptfont\bffam=\sevenbf
    \textfont\ttfam=\fourteentt
    \textfont\cpfam=\fourteencp
    \textfont\ssfam=\fourteenss     \scriptfont\ssfam=\tenss
      \scriptscriptfont\ssfam=\sevenss
    \samef@nt
    \b@gheight=14pt
    \setbox\strutbox=\hbox{\vrule height 0.85\b@gheight
                                depth 0.35\b@gheight width\z@ }}
\def\twelvepoint{\relax
    \textfont0=\twelverm          \scriptfont0=\ninerm
      \scriptscriptfont0=\sixrm
    \textfont1=\twelvei           \scriptfont1=\ninei
      \scriptscriptfont1=\sixi
    \textfont2=\twelvesy           \scriptfont2=\ninesy
      \scriptscriptfont2=\sixsy
    \textfont3=\twelveex          \scriptfont3=\ninex
      \scriptscriptfont3=\sixex
    \textfont\itfam=\twelveit    %\scriptfont\itfam=\nineit
    \textfont\slfam=\twelvesl    %\scriptfont\slfam=\ninesl
      \scriptscriptfont\slfam=\sixsl
    \textfont\bffam=\twelvebf     \scriptfont\bffam=\ninebf
      \scriptscriptfont\bffam=\sixbf
    \textfont\ttfam=\twelvett
    \textfont\cpfam=\twelvecp
    \textfont\ssfam=\twelvess     \scriptfont\ssfam=\niness
      \scriptscriptfont\ssfam=\sixss
    \samef@nt
    \b@gheight=12pt
    \setbox\strutbox=\hbox{\vrule height 0.85\b@gheight
                                depth 0.35\b@gheight width\z@ }}
\def\tenpoint{\relax
    \textfont0=\tenrm          \scriptfont0=\sevenrm
      \scriptscriptfont0=\fiverm
    \textfont1=\teni           \scriptfont1=\seveni
      \scriptscriptfont1=\fivei
    \textfont2=\tensy          \scriptfont2=\sevensy
      \scriptscriptfont2=\fivesy
    \textfont3=\tenex          \scriptfont3=\sevenex
      \scriptscriptfont3=\fivex
    \textfont\itfam=\tenit     \scriptfont\itfam=\seveni
    \textfont\slfam=\tensl     \scriptfont\slfam=\sevensl
      \scriptscriptfont\slfam=\fivesl
    \textfont\bffam=\tenbf     \scriptfont\bffam=\sevenbf
      \scriptscriptfont\bffam=\fivebf
    \textfont\ttfam=\tentt
    \textfont\cpfam=\tencp
    \textfont\ssfam=\tenss     \scriptfont\ssfam=\sevenss
      \scriptscriptfont\ssfam=\fivess
    \samef@nt
    \b@gheight=10pt
    \setbox\strutbox=\hbox{\vrule height 0.85\b@gheight
                                depth 0.35\b@gheight width\z@ }}
%
%%%%%%%%%%%%%%%%%%%%%%%%%%%%%%%%%%%%%%%%%%%%%%%%%%%%%%%%%%%%%%%%%%%%%%%%
%
%   Next, I define basic spacing parameters.
%
\normalbaselineskip = 15pt plus 0.2pt minus 0.1pt %was 20pt ...
\normallineskip = 1.5pt plus 0.1pt minus 0.1pt
\normallineskiplimit = 1.5pt
\newskip\normaldisplayskip
\normaldisplayskip = 15pt plus 5pt minus 10pt %was 20pt ...
\newskip\normaldispshortskip
\normaldispshortskip = 6pt plus 5pt
\newskip\normalparskip
\normalparskip = 6pt plus 2pt minus 1pt
\newskip\skipregister
\skipregister = 5pt plus 2pt minus 1.5pt
\newif\ifsingl@    \newif\ifdoubl@
\newif\iftwelv@    \twelv@true
\def\singlespace{\singl@true\doubl@false\spaces@t}
\def\doublespace{\singl@false\doubl@true\spaces@t}
\def\normalspace{\singl@false\doubl@false\spaces@t}
\def\Tenpoint{\tenpoint\twelv@false\spaces@t}
\def\Twelvepoint{\twelvepoint\twelv@true\spaces@t}
\def\spaces@t{\relax
      \iftwelv@ \ifsingl@\subspaces@t3:4;\else\subspaces@t1:1;\fi
       \else \ifsingl@\subspaces@t3:5;\else\subspaces@t4:5;\fi \fi
      \ifdoubl@ \multiply\baselineskip by 5
         \divide\baselineskip by 4 \fi }
\def\subspaces@t#1:#2;{
      \baselineskip = \normalbaselineskip
      \multiply\baselineskip by #1 \divide\baselineskip by #2
      \lineskip = \normallineskip
      \multiply\lineskip by #1 \divide\lineskip by #2
      \lineskiplimit = \normallineskiplimit
      \multiply\lineskiplimit by #1 \divide\lineskiplimit by #2
      \parskip = \normalparskip
      \multiply\parskip by #1 \divide\parskip by #2
      \abovedisplayskip = \normaldisplayskip
      \multiply\abovedisplayskip by #1 \divide\abovedisplayskip by #2
      \belowdisplayskip = \abovedisplayskip
      \abovedisplayshortskip = \normaldispshortskip
      \multiply\abovedisplayshortskip by #1
        \divide\abovedisplayshortskip by #2
      \belowdisplayshortskip = \abovedisplayshortskip
      \advance\belowdisplayshortskip by \belowdisplayskip
      \divide\belowdisplayshortskip by 2
      \smallskipamount = \skipregister
      \multiply\smallskipamount by #1 \divide\smallskipamount by #2
      \medskipamount = \smallskipamount \multiply\medskipamount by 2
      \bigskipamount = \smallskipamount \multiply\bigskipamount by 4 }
\def\normalbaselines{ \baselineskip=\normalbaselineskip
   \lineskip=\normallineskip \lineskiplimit=\normallineskip
   \iftwelv@\else \multiply\baselineskip by 4 \divide\baselineskip by 5
     \multiply\lineskiplimit by 4 \divide\lineskiplimit by 5
     \multiply\lineskip by 4 \divide\lineskip by 5 \fi }
\Twelvepoint  % That's the default
%%%%%%%%%%%%%%%%%%%%%%%%%%%%\Tenpoint   % Not in Leuven it wasn't!
%
\interlinepenalty=50
\interfootnotelinepenalty=5000
\predisplaypenalty=9000
\postdisplaypenalty=500
\hfuzz=1pt
\vfuzz=0.2pt
\dimen\footins=24 truecm % 8 truein in SB
\hoffset=0 truemm % 0 in SB, 6.5mm in Leuven
\voffset=0 truemm % 0in in SB, 5 truemm in Leuven
%
%%%%%%%%%%%%%%%%%%%%%%%%%%%%%%%%%%%%%%%%%%%%%%%%%%%%%%%%%%%%%%%%%%
%
% Now some output macros
%
%%%%%%%%%%%%%%%%%%%%%%%%%%%%%%%%%%%%%%%%%%%
%
% footnote numbering macros
%
%
\def\footnote#1{\edef\@sf{\spacefactor\the\spacefactor}#1\@sf
      \insert\footins\bgroup\singl@true\doubl@false\Tenpoint
      \interlinepenalty=\interfootnotelinepenalty \let\par=\endgraf
        \leftskip=\z@skip \rightskip=\z@skip
        \splittopskip=10pt plus 1pt minus 1pt \floatingpenalty=20000
        \smallskip\item{#1}\bgroup\strut\aftergroup\@foot\let\next}
\skip\footins=\bigskipamount % space added when footnote is present
\dimen\footins=24truecm % maximum footnotes per page (8 truein in USA)
\newcount\fnotenumber
\def\clearfnotenumber{\fnotenumber=0}
\def\fnote{\advance\fnotenumber by1 \footnote{$^{\the\fnotenumber}$}}
\clearfnotenumber
%
% section and appendix macros
%
\newcount\secnumber
\newcount\appnumber
\newif\ifs@c % this is true if within a section as opposed to an appendix
\newif\ifs@cd % this is true if the article is being section'd
\s@cdtrue % this is the default
\def\unsectioned{\s@cdfalse\let\section=\subsection}
\def\clearappnumber{\appnumber=64}
\def\clearsecnumber{\secnumber=0}
\newskip\sectionskip         \sectionskip=\medskipamount
\newskip\headskip            \headskip=8pt plus 3pt minus 3pt
\newdimen\sectionminspace    \sectionminspace=10pc
\newdimen\referenceminspace  \referenceminspace=25pc
\def\Titlestyle#1{\par\begingroup \interlinepenalty=9999
     \leftskip=0.02\hsize plus 0.23\hsize minus 0.02\hsize
     \rightskip=\leftskip \parfillskip=0pt
     \advance\baselineskip by 0.5\baselineskip%this is a test...
     \hyphenpenalty=9000 \exhyphenpenalty=9000
     \tolerance=9999 \pretolerance=9000
     \spaceskip=0.333em \xspaceskip=0.5em
     \seventeenpoint
  \noindent #1\par\endgroup }
\def\titlestyle#1{\par\begingroup \interlinepenalty=9999
     \leftskip=0.02\hsize plus 0.23\hsize minus 0.02\hsize
     \rightskip=\leftskip \parfillskip=0pt
     \hyphenpenalty=9000 \exhyphenpenalty=9000
     \tolerance=9999 \pretolerance=9000
     \spaceskip=0.333em \xspaceskip=0.5em
     \fourteenpoint
   \noindent #1\par\endgroup }%  the \Npoint only takes care of spacing.
%                                a font is always specified when calling this
%                                macro.  In computers with little room for
%                                character-size data it is convenient to % out
%                                all the font definitions from the \Npoint
%                                macros.
%
\def\spacecheck#1{\dimen@=\pagegoal\advance\dimen@ by -\pagetotal
   \ifdim\dimen@<#1 \ifdim\dimen@>0pt \vfil\break \fi\fi}
\def\section#1{\cleareqnumber \s@ctrue \global\advance\secnumber by1
   \message{Section \the\secnumber: #1}
   \par \ifnum\the\lastpenalty=30000\else
   \penalty-200\vskip\sectionskip \spacecheck\sectionminspace\fi
   \noindent {\caps\enspace\S\the\secnumber\quad #1}\par
   \nobreak\vskip\headskip \penalty 30000 }
\def\subsection#1{\par
   \ifnum\the\lastpenalty=30000\else \penalty-100\smallskip
   \spacecheck\sectionminspace\fi
   \noindent\undertext{#1}\enspace \vadjust{\penalty5000}}

\def\undertext#1{\vtop{\hbox{#1}\kern 1pt \hrule}}
\def\subsubsection#1{\par
   \ifnum\the\lastpenalty=30000\else \penalty-100\smallskip \fi
   \noindent\hbox{#1}\enspace \vadjust{\penalty5000}}

\def\appendix#1{\cleareqnumber \s@cfalse \global\advance\appnumber by1
   \message{Appendix \char\the\appnumber: #1}
   \par \ifnum\the\lastpenalty=30000\else
   \penalty-200\vskip\sectionskip \spacecheck\sectionminspace\fi
   \noindent {\caps\enspace Appendix \char\the\appnumber\quad #1}\par
   \nobreak\vskip\headskip \penalty 30000 }
\clearsecnumber
\clearappnumber
%
% macros for references, acknowledgements, and note added
%
\def\ack{\par\penalty-100\medskip \spacecheck\sectionminspace
   \line{\iftwelv@\fourteencp\else\twelvecp\fi\hfil {\caps
Acknowledgements}\hfil}%
\nobreak\vskip\headskip }
\def\refs{\begingroup \par\penalty-100\medskip \spacecheck\sectionminspace
   \line{\iftwelv@\fourteencp\else\twelvecp\fi\hfil REFERENCES\hfil}%
\nobreak\vskip\headskip \frenchspacing }
\def\endrefs{\par\endgroup}
%--- Note added
%
\newcount\refnumber
\def\clearrefnumber{\refnumber=0}  \clearrefnumber
\newwrite\R@fs                              %This opens a file .refs with
\immediate\openout\R@fs=\jobname.references %the references in order of
                                            %appearance.
\def\closerefs{\immediate\closeout\R@fs} %close file so that TeX can read it
\def\refsout{\closerefs\refs
\catcode`\@=11                          % we must do this since the
\input\jobname.references               % references expand to
\catcode`\@=12			        % primitives containing @'s
\endrefs}
\def\refitem#1{\item{{\bf #1}}}%just bolds it so that \bf does not expand
\def\ifundefined#1{\expandafter\ifx\csname#1\endcsname\relax}
%
%  new reference macros.  Now just say \[_label], and in the
%  ref_defs.tex file type \refdef[_label]{reference}
%
\def\[#1]{\ifundefined{#1R@FNO}%
\global\advance\refnumber by1%
\expandafter\xdef\csname#1R@FNO\endcsname{[\the\refnumber]}%
\immediate\write\R@fs{\noexpand\refitem{\csname#1R@FNO\endcsname}%
\noexpand\csname#1R@F\endcsname}\fi{\bf \csname#1R@FNO\endcsname}}
\def\refdef[#1]#2{\expandafter\gdef\csname#1R@F\endcsname{{#2}}}
%
% equation numbering macros
%
%  better than before.  just do \(_label) both to refer or to
%                         define.
%
%        the generic equation is \()
%
%
%
\newcount\eqnumber
\def\cleareqnumber{\eqnumber=0}
\newif\ifal@gn \al@gnfalse  % this is true if within an \eqalignno
% at some point try the following:
%\def\eqnalign#1{\al@gntrue \vbox{\eqalignno{#1}} \al@gnfalse}
% but meanwhile let`s define a new macro...
\def\veqnalign#1{\al@gntrue \vbox{\eqalignno{#1}} \al@gnfalse}
\def\eqnalign#1{\al@gntrue \eqalignno{#1} \al@gnfalse}
\def\(#1){\relax%
\ifundefined{#1@Q}
 \global\advance\eqnumber by1
 \ifs@cd
  \ifs@c
   \expandafter\xdef\csname#1@Q\endcsname{{%
\noexpand\rm(\the\secnumber .\the\eqnumber)}}
  \else
   \expandafter\xdef\csname#1@Q\endcsname{{%
\noexpand\rm(\char\the\appnumber .\the\eqnumber)}}
  \fi
 \else
  \expandafter\xdef\csname#1@Q\endcsname{{\noexpand\rm(\the\eqnumber)}}
 \fi
 \ifal@gn
    & \csname#1@Q\endcsname
 \else
    \eqno \csname#1@Q\endcsname
 \fi
\else%
\csname#1@Q\endcsname\fi\global\let\@Q=\relax}
%
% macros for running heads and page numbering
%
\newif\iffrontpage \frontpagefalse
\headline={\hfil}
\footline={\iffrontpage\hfil\else \hss\twelverm
-- \folio\ --\hss \fi }
\def\monthname{\relax\ifcase\month 0/\or January\or February\or
   March\or April\or May\or June\or July\or August\or September\or
   October\or November\or December\else\number\month/\fi}
\hsize=16 truecm
\vsize=22 truecm
\skip\footins=\bigskipamount
\normalspace
%
%%%%%%%%%%%%%%%%%%%%%%%%%%%%%%%%%%%%%%%%%%%%%%%%%%%%%%%%%%%%%%%%%%%%%%%
%
%   Here come macros for title pages.
%
\newskip\frontpageskip
\newif\ifp@bblock \p@bblocktrue
\newif\ifm@nth \m@nthtrue
\newtoks\pubnum
\newtoks\pubtype
\newtoks\m@nthn@me
\newcount\Ye@r
\advance\Ye@r by \year
\advance\Ye@r by -1900
\def\Year#1{\Ye@r=#1}%--- set the year by hand
\def\Month#1{\m@nthfalse \m@nthn@me={#1}}
\def\m@nthname{\ifm@nth\monthname\else\the\m@nthn@me\fi}
\def\titlepage{\global\frontpagetrue\hrule height\z@ \relax
               \ifp@bblock\pubblock\fi\relax }
\def\endtitlepage{\vfil\break
                  \frontpagefalse} %I took a \pageno=1 from here
\def\bonntitlepage{\global\frontpagetrue\hrule height\z@ \relax
               \ifp@bblock\pubblock\fi\relax }
\frontpageskip=12pt plus .5fil minus 2pt
\pubtype={\iftwelv@\twelvesl\else\tensl\fi\ (Preliminary Version)}
\pubnum={?}
\def\nopubblock{\p@bblockfalse}
\def\pubblock{\line{\hfil\iftwelv@\twelverm\else\tenrm\fi%
US--\number\Ye@r--\the\pubnum\the\pubtype}
              \line{\hfil\iftwelv@\twelverm\else\tenrm\fi%
\m@nthname\ \number\year}}
\def\title#1{\vskip\frontpageskip\Titlestyle{\caps #1}\vskip3\headskip}
%                                 ^---notice capital Titlestyle...
\def\author#1{\vskip.5\frontpageskip\titlestyle{\caps #1}\nobreak}
\def\andauthor{\vskip.5\frontpageskip\centerline{and}\author}

\def\address#1{\par\kern 5pt\titlestyle{%\iftwelv@\twelvepoint\else\tenpoint\fi
\it #1}}
\def\andaddress{\par\kern 5pt \centerline{\sl and} \address}

\def\Santiago{\address{Departamento de F{\'\i}sica de Part{\'\i}culas
Elementales\break
Universidad de Santiago, Santiago de Compostela, 15706, SPAIN}}

\def\abstract#1{\par\dimen@=\prevdepth \hrule height\z@ \prevdepth=\dimen@
   \vskip\frontpageskip\spacecheck\sectionminspace
   \centerline{\iftwelv@\fourteencp\else\twelvecp\fi ABSTRACT}\vskip\headskip
   {\noindent #1}}
%

%
%%%%%%%%%%%%%%%%%%%%%%%%%%%%%%%%%%%%%%%%%%%%%%%%%%%%%%%%%%%%%%%%%%
%
% macros for leaders, boxes, underline, ...
%
\def\leaderfill{\leaders\hbox to 1em{\hss.\hss}\hfill}%--- leading ...
 %--- underline
\def\boxit#1{\vcenter{\hrule\hbox{\vrule\kern8pt
      \vbox{\kern8pt#1\kern8pt}\kern8pt\vrule}\hrule}}%--- box
 %--- box in $$...$$

%
%%%%%%%%%%%%%%%%%%%%%%%%%%%%%%%%%%%%%%%%%%%%%%%%%%%%%%%%%%%%%%%%%%
%
%  Now come basic non-math macros
%
\def\ref#1{{\bf [#1]}}%--- [ref]
%--- et al.
\def\ie{{\it i.e.\/}}%--- i.e.
%--- e.g.
%--- Cf.
%--- cf.
 %--- double left quote
%--- th as in fifth
\def\nl{\hfil\break}%--- new line
%--- just an abbrev.
%--- 1/2
%
%%%%%%%%%%%%%%%%%%%%%%%%%%%%%%%%%%%%%%%%%%%%%%%%%%%%%%%%%%%%%%%%
%
% Now some math macros
%
%%%%%%%%%%%%%%%%%%%%%%%%%%%%%%%%%%%%%
%
% First macros for theorems, definitions, ...
%
\newif\ifm@thstyle \m@thstylefalse
\def\mathstyle{\m@thstyletrue}
\def\proclaim#1#2\par{\smallbreak\begingroup%        small --> med???
\advance\baselineskip by -0.25\baselineskip%
\advance\belowdisplayskip by -0.35\belowdisplayskip%
\advance\abovedisplayskip by -0.35\abovedisplayskip%
    \noindent{\caps#1.\enspace}{#2}\par\endgroup%
\smallbreak}%--- defs, thms, ...                     small --> med???
\def\m@kem@th<#1>#2#3{%
\ifm@thstyle \global\advance\eqnumber by1
 \ifs@cd
  \ifs@c
   \expandafter\xdef\csname#1\endcsname{{%
\noexpand #2\ \the\secnumber .\the\eqnumber}}
  \else
   \expandafter\xdef\csname#1\endcsname{{%
\noexpand #2\ \char\the\appnumber .\the\eqnumber}}
  \fi
 \else
  \expandafter\xdef\csname#1\endcsname{{\noexpand #2\ \the\eqnumber}}
 \fi
 \proclaim{\csname#1\endcsname}{#3}
\else
 \proclaim{#2}{#3}
\fi}
%
%
%  To use the new Math macros...
%
%         \Thm<_label>{Statemenet of the Thm} etc...
%
%     where _label is the label of the Thm.  The generic
%     Thm has an empty label.
%
%     To refer to it, just say \<_label>
%
\def\Thm<#1>#2{\m@kem@th<#1M@TH>{Theorem}{\sl#2}}%--- Theorem
\def\Prop<#1>#2{\m@kem@th<#1M@TH>{Proposition}{\sl#2}}%--- Proposition
\def\Def<#1>#2{\m@kem@th<#1M@TH>{Definition}{\rm#2}}%--- Definition
\def\Lem<#1>#2{\m@kem@th<#1M@TH>{Lemma}{\sl#2}}%--- Lemma
\def\Cor<#1>#2{\m@kem@th<#1M@TH>{Corollary}{\sl#2}}%--- Corollary
\def\Conj<#1>#2{\m@kem@th<#1M@TH>{Conjecture}{\sl#2}}%--- Conjecture
\def\Rmk<#1>#2{\m@kem@th<#1M@TH>{Remark}{\rm#2}}%--- Remark
\def\Exm<#1>#2{\m@kem@th<#1M@TH>{Example}{\rm#2}}%--- Example
\def\Qry<#1>#2{\m@kem@th<#1M@TH>{Query}{\it#2}}%--- Query
\def\Example<#1>#2{\m@kem@th<#1M@TH>{Example}{\it#2}}%--- Example
\def\Exercise<#1>#2{\m@kem@th<#1M@TH>{Execise}{\it#2}}%--- Exercise
%
%--- Proof
%
\let\Pf=\Proof
\let\Example=\Exm
\def\<#1>{\csname#1M@TH\endcsname}
%
% We then continue with basic mathematics
%
%--- def over =
\def\qed{\vrule width 0.6em height 0.5em depth 0.2em}%--- Halmos Q.E.D.
\def\QED{\enspace\qed}
%--- implies
%--- is implied by
%--- if and only if
\def\lapprox{\hbox{\lower3pt\hbox{$\buildrel<\over\sim$}}}% approx lt
\def\gapprox{\hbox{\lower3pt\hbox{$\buildrel<\over\sim$}}}% approx gt
%--- exponential
\def\quotient#1#2{#1/\lower0pt\hbox{${#2}$}}%--- factor objects
%
% Arrow stuff
%
%--- injective map
%--- surjective map
%--- bijective map
\def\to{\rightarrow}%--- mapping
%--- long mapping
%--- isom over -->
%--- just an abbrev.
%

%
 %--- commutative diagram macro
%-- lin map over arrow
 %--- map in complex
%
% Numbers...
%
 %--- reals
\def\comps{{\bf C}} %--- complex nos.
 %--- quaternions
\def\integ{{\bf Z}} %--- integers
 %--- rationals
 %--- naturals
 %--- ground field
%
% Algebra
%
%--- Hom(omorphisms)
%--- tr(ace)
\def\Tr{{\rm Tr}~}%--- Tr(ace)
%--- End(omorphisms)
%--- Mor(phisms)
%--- Aut(omorphisms)
%--- aut(omorphisms)
%--- supertrace
%--- superdeterminant
%--- kernel
%--- cokernel
%--- image
\def\underrightarrow#1{\vtop{\ialign{##\crcr
      $\hfil\displaystyle{#1}\hfil$\crcr
      \noalign{\kern-\p@\nointerlineskip}
      \rightarrowfill\crcr}}} %--- modification of \overrightarrow
\def\underleftarrow#1{\vtop{\ialign{##\crcr
      $\hfil\displaystyle{#1}\hfil$\crcr
      \noalign{\kern-\p@\nointerlineskip}
      \leftarrowfill\crcr}}}  %--- modification of \overleftarrow

%
% Brackets,...
%
\def\comm#1#2{\left[#1\, ,\,#2\right]}%--- [ , ]
%--- { , }
\def\pbr#1#2{\left\{#1\, ,\,#2\right\}}       %--- { , }
%--- [ , }
%--- structure const.
%
% Analysis anyone?
%
%--- Lie derivative
% lft var derivative
% rgt var derivative
%--- vartnl derivative
% double vartl deriv.
%--- rgt prtl derivative
%--- partial derivative
%
%                                                           dble partl deriv.
%--- full lft derivative
%--- full rgt derivative
\def\der#1#2{{{d #1}\over {d #2}}}%--- full derivative
%--- laplacian
%--- cov. ext. der.
%
% Dirac slashes
%
%--- D slash
%--- del slash
%--- A slash
\def\wt{\widetilde} 
%
%%%%%%%%%%%%%%%%%%%%%%%%%%%%%%%%%%%%%%%%%%%%%%%%%%%%%%%%%%%%%%%%
%
%  These are the macros to make Young tableaux
%
\newdimen\unit
\newdimen\redunit
%
%   this puts the ref. point of #1 at coordinates (#2,#3)
%
\def\p@int#1:#2 #3 {\rlap{\kern#2\unit
     \raise#3\unit\hbox{#1}}}
%
% this defines the sides of the tableau
% notice that \rver and \lver coincide
% It would have been natural for \rver to have negative
% width but that does not print in TeX...
%
\def\th@r{\vrule height0\unit depth.1\unit width1\unit}
\def\bh@r{\vrule height.1\unit depth0\unit width1\unit}
\def\lv@r{\vrule height1\unit depth0\unit width.1\unit}
\def\rv@r{\vrule height1\unit depth0\unit width.1\unit}
%
% this is the tableau: the .9 is due to the unnatural definition
% of \rver
%
\def\t@ble@u{\hbox{\p@int\bh@r:0 0
                   \p@int\lv@r:0 0
                   \p@int\rv@r:.9 0
                   \p@int\th@r:0 1
                   }
             }
%
% we now define the tableau at a particular location
%
\def\t@bleau#1#2{\rlap{\kern#1\redunit
     \raise#2\redunit\t@ble@u}}
%
%  Now a macro to make a column of #1 tableaux down at (#2,#3)
%
\newcount\n
\newcount\m
\def\makecol#1#2#3{\n=0 \m=#3
  \loop\ifnum\n<#1{}\advance\m by -1 \t@bleau{#2}{\number\m}\advance\n by 1
\repeat}
%
%   Now a macro to make a row of #1 tableaux at (#2,#3) to the right
%
\def\makerow#1#2#3{\n=0 \m=#3
 \loop\ifnum\n<#1{}\advance\m by 1 \t@bleau{\number\m}{#2}\advance\n by 1
\repeat}
%
% Some useful ready made Young tableaux
%
\def\checkunits{\ifinner \unit=6pt \else \unit=8pt \fi
                \redunit=0.9\unit } %these are the basic sizes
\def\ytsym#1{\checkunits\kern-.5\unit
  \vcenter{\hbox{\makerow{#1}{0}{0}\kern#1\unit}}\kern.5em} % #1 symmetrized
%                                                             tableaux
\def\ytant#1{\checkunits\kern.5em
  \vcenter{\hbox{\makecol{#1}{0}{0}\kern1\unit}}\kern.5em} % #1 antisymmetrized
%                                                            tableaux
\def\ytwo#1#2{\checkunits
  \vcenter{\hbox{\makecol{#1}{0}{0}\makecol{#2}{1}{0}\kern2\unit}}
                  \ } % 2 column #1 #2 (left->right) Young tableau
\def\ythree#1#2#3{\checkunits
  \vcenter{\hbox{\makecol{#1}{0}{0}\makecol{#2}{1}{0}\makecol{#3}{2}{0}%
\kern3\unit}}
                  \ } % 3 column #1 #2 #3 (left->right) Young tableau

%%%%%%%%%%%%%%%%%%%%%%%%%%%%%%%%%%%%%%%%%%%%%%%%%%%%%%%%%%%%%%%%%%%%%
%
% Finally some useful macros for journals, ...
%
%
\def\PRL#1#2#3{{\sl Phys. Rev. Lett.} {\bf#1} (#2) #3}
\def\NPB#1#2#3{{\sl Nucl. Phys.} {\bf B#1} (#2) #3}

\def\CMP#1#2#3{{\sl Comm. Math. Phys.} {\bf #1} (#2) #3}

\def\PLB#1#2#3{{\sl Phys. Lett.} {\bf #1B} (#2) #3}
\def\JMP#1#2#3{{\sl J. Math. Phys.} {\bf #1} (#2) #3}

\def\PR#1#2#3{{\sl Phys. Reports} {\bf #1} (#2) #3}

\def\FAaIA#1#2#3{{\sl Funct. Anal. Appl.} {\bf #1} (#2)
#3}

\def\Invm#1#2#3{{\sl Invent. math.} {\bf #1} (#2) #3}
\def\LMP#1#2#3{{\sl Letters in Math. Phys.} {\bf #1} (#2) #3}
\def\IJMPA#1#2#3{{\sl Int. J. Mod. Phys.} {\bf A#1} (#2) #3}

\def\TMP#1#2#3{{\sl Theor. Mat. Phys.} {\bf #1} (#2) #3}

\def\JSM#1#2#3{{\sl J. Soviet Math.} {\bf #1} (#2) #3}
\def\MPLA#1#2#3{{\sl Mod. Phys. Lett.} {\bf A#1} (#2) #3}

\def\PJAS#1#2#3{{\sl Proc. Jpn. Acad. Sci.} {\bf #1} (#2) #3}
\def\JPSJ#1#2#3{{\sl J. Phys. Soc. Jpn.} {\bf #1} (#2) #3}
\def\JETPL#1#2#3{{\sl  Sov. Phys. JETP Lett.} {\bf #1} (#2) #3}

\catcode`\@=12 % @ no longer a letter
%
%
%
%%%%%%%%%%%%%%%%%%%%%%%%%%%%%%%%%%%%%%%%%%%%%%%
%   These are the local macros for q-W
%%%%%%%%%%%%%%%%%%%%%%%%%%%%%%%%%%%%%%%%%%%%%%%

\def\qWKP{ q \hbox{-}{\ssf W}_{\rm KP}}
\def\WKP{{\ssf W}_{\rm KP}}
\def\max{{\rm max}}
\def\min{{\rm min}}

\def\inv{^{-1}}
\def\G{{\cal G}}
\def\F{F}
\def\g{{\mit g}}
\def\M{{\cal M}}
\def\R{{\cal R}}
\def\J{ J }
\def\P{P}
\def\p{{\ssf p}}
\def\lbc{\lbrace}
\def\rbc{\rbrace}

% kth partial

\def\d{\partial}
\def\pdo{{\hbox{$\Psi$DO}}}

\def\comb[#1/#2]{\left[{#1\atop#2}\right]}
\def\pair#1#2{\langle #1,#2\rangle} %--- dual pairing

\def\rflecha#1{\setbox1=\hbox{\ninerm ~#1~}%
\buildrel\hbox{\ninerm #1}\over{\hbox to\wd1{\rightarrowfill}}}
\def\lflecha#1{\setbox1=\hbox{\ninerm ~#1~}%
\buildrel\hbox{\ninerm #1}\over{\hbox to\wd1{\leftarrowfill}}}

\def\res{{\rm \,res}}

\def\ResT{\res_{T}\,}
\def\ResD{\res_{\d_q}\,}
\def\comb[#1/#2]{\left[{#1\atop#2}\right]}
\def\qcomb[#1/#2]{\left[{#1\atop#2}\right]_{\! q}}

\def\ov{\over}
\def\fr#1/#2{\mathord{\hbox{${#1}\over{#2}$}}}

\def\med{\fr1/2}

\def\W{\mathord{\ssf W}}

\def\ope[#1][#2]{{{#2}\over{\ifnum#1=1 {z-w} \else {(z-w)^{#1}}\fi}}}

%--- EM tensor for ghost system

\def\suma2#1#2{\sum_{#1=0}^{#2}}

%%%%%%%%%%%%%%%%%%%%%%%%%%%%%%%%%%%%%%%%%%%%%%%%%%%%%%%%%%
% These are the local refs for L=AB
%%%%%%%%%%%%%%%%%%%%%%%%%%%%%%%%%%%%%%%%%%%%%%%%%%%%%%%%%%

\refdef[WReview]{P.~Bouwknegt and K.~Schoutens, {\sl ${\cal
W}$-Symmetry in Conformal Field Theory},  {\it Phys. Reps.} to
appear.}
\refdef[Univ]{J.~M.~Figueroa-O'Farrill and E.~Ramos,
\JMP{33}{1992}{833}.}
\refdef[Wn]{A.~B.~Zamolodchikov, \TMP{65}{1986}{1205};\nl
V.~A.~Fateev and S.~L.~Lykyanov, \IJMPA{3}{1988}{507}.}
\refdef[Dickey]{L.~A.~Dickey,  {\sl Soliton equations and Hamiltonian
systems},  Advanced Series in Mathematical Physics Vol.12,  World
Scientific Publ.~Co..}
\refdef[GD]{I.~M.~Gel'fand and L.~A.~Dickey, {\sl A family of
Hamiltonian structures connected with integrable nonlinear
differential equations}, Preprint 136, IPM AN SSSR, Moscow (1978).}
\refdef[Adler]{M.~Adler, \Invm{50}{1979}{403}.}
\refdef[KP]{E.~Date, M.~Jimbo, M.~Kashiwara, and T.~Miwa
\PJAS{57A}{1981}{387}; \JPSJ{50}{1981}{3866}.}
\refdef[WKP]{L.~A.~Dickey, {\sl Annals NY Acad.~Sci.} {\bf 491}(1987)
131;\nl J.~M.~Figueroa-O'Farrill, J.~Mas, and E.~Ramos,
\PLB{266}{1991}{298};\nl F.~Yu and Y.-S.~Wu, \NPB{373}{1992}{713}.}
\refdef[WKPq]{J.~M.~Figueroa-O'Farrill, J.~Mas, and E.~Ramos
\CMP{158}{1993}{17}.}
\refdef[WinftyKP]{K.~Yamagishi, \PLB{259}{1991}{436};\nl F.~Yu and
Y.-S.~Wu, \PLB{236}{1991}{220}}
\refdef[Schoutens]{P. Bouwknegt and K. Schoutens,~{\it
$W$ symmetry in conformal field theory},~\PR{223}{1993}{183-276}.}
\refdef[Wgeom]{J.~Gervais and Y.~ Matsuo,\PLB{282}{1992}{309},~({\tt
hep-th/9110028}).}
\refdef[Class]{J.~M.~Figueroa-O'Farrill and E.~Ramos, {\it
The classical limit of $\W$ algebras},~\PLB{282}{1992}{357},~ ({\tt
hep-th/9202040}).}
\refdef[FigRaSo]{ J.~M.~Figueroa-O'Farrill, E.~Ramos and S.~ Stanciu,
{\it A geometrical interpretation of classical W transformations},
\PLB{297}{1992}{289}}
\refdef[RaRo]{ E.~Ramos and J.~Roca, {\it Extended gauge
invariance in geometrical particle models and the geometry of W
symmetry}. \NPB{452}{1995}{705}.}
\refdef[Radul]{A.~O.~Radul, in {\sl Applied methods of nonlinear
analysis and control}, pp. 149-157, Mironov, Moroz, and Tshernjatin,
eds.,  MGU 1987 (Russian).}
\refdef[winfty]{I.~Bakas, \PLB{228}{1989}{406}; \CMP{134}{1990}{487}.}
\refdef[Winfty]{C.~N.~Pope, L.~J.~Romans, and X.~Shen,
\NPB{339}{1990}{191}.}
\refdef[Woneplusinfty]{C. N. Pope, L. J. Romans, and X. Shen,
\PLB{242}{1990}{401}.}
\refdef[RadulFigRa]{A.O.~ Radul, \JETPL{50}{1989}{371},\nl
J.M.~ Figueroa O'Farrill and E. Ramos,~ \LMP{27}{1993}{223}~({\tt
hep-th/9211036}).
}
\refdef[KeZa]{B.A.~Khesin adn I.S.~ Zakharevich,{\it Poisson-Lie group od
pseudodifferential operators and fractional KP-KdV hierarchies},~{\sl C.R.
Acad. Sci.,
{\bf 315} S\'er.I (1993), 621-626},\nl
{\it Poisson-Lie group of pseudodifferential operators},~
\CMP{171}{1995}{475-530},~ ({\tt hep-th/9312088}).}
\refdef[DS]{V.~G.~Drinfel'd and V.~V.~Sokolov, \JSM{30}{1984}{1975}.}
\refdef[redukp]{
W. Oevel and W. Strampp,~  \CMP{157}{1993}{51-81}.\nl
F.~Yu and Y.-S.~Wu, \PRL{68}{1992}{2996} ({\tt
hep-th/9112009}); \nl
F. Yu, \LMP{29}{1993}{175};\nl
H.~Aratyn, L.~A.~Ferreira, J.~F.~Gomes, and
A.~H.~Zimerman, ~\NPB{402}{1993}{85-117} \nl
L.~Bonora, Q.P.~Liu and C.S.~ Xiong,~ {\sl The
integrable hierarchy constructed from a pair of KdV type
hierarchies and its
associated $\W$ algebra}, ( {\tt hep-th/9408035}); \nl
D. Depireux and J. Shiff, preprint  IASSNS-HEP-92-66. }
\refdef[dickeyred]{L.~A.~Dickey,
{\sl On the constrained KP hierarchy, I and II}, ( {\tt hep-th/9407038} and
{\tt hep-th/9411005}).}
\refdef[hidi]{F. Mart\'\i nez-Mor\'as and E. Ramos,~
 \CMP{157}{1993}{573-589}.
\nl
F. Mart\'\i nez-Mor\'as, J. Mas and E. Ramos,
{}~\MPLA{8}{1993}{2189-2197}.}
\refdef[OeRa]{W.~Oevel and O.~Ragnisco,~ {\sl Physica }{\bf A} 161 (1989)
181-220.}
\refdef[KLR]{B.~Khesin, V.~Lyubashenko and C.~Roger, {\it Extensions and
contractions
of the Lie algebra of $q$-pseudodifferential symbols}, ~{\tt hep-th/9403189}.}
\refdef[Frenkel]{E. ~Frenkel, {\it Deformations of the KdV hierarchy and
related
soliton equations} ~{\tt q-alg/9511003} .}
\refdef[QuantW]{
J. ~Shiraishi,~H.~Kubo,~H.~ Awata,S.~Odake,{\it A quantum deformation of the
Virasoro algebra and Macdonald symmetric functions},~{\tt q/alg-9507034}.\nl
S. Lykyanov and Ya. Pugay,~{\it Bosonization of ZF algebras: direction toward
deformed Virasoro algebra},~{\tt hep-th/9412128}.\nl
R. Kemmoku and S. Saito, {\it $W_{1+\infty}$ as a discretization of the
Virasoro
algebra},~{\tt hep-th/9411027}. \nl
H-T.~Sato,~{\it $q$-Virasoro Operators from $q$-Noether Currents},~
{\tt hep-th/9510189}.
}
\refdef[FResh]{E. Frenkel and N. Reshetikhin,{\it  Quantum affine
algebras and  deformations of the Virasoro and W-algebras},~{\tt
q-alg/9505025}. }
\refdef[Babelon]{O.~ Babelon {\it Exchange formula and lattice deformation of
the
Virasoro algebra},~\PLB{238}{1990}{234-238}.}
\refdef[kupwilson]{B.A. ~Kupershmidt and G.~ Wilson,~ \Invm{62}{1981}{403}.}
\refdef[Hall]{A. Capelli, C. Trugenberger and G. Zemba, \NPB{396}{1993}{465}
\nl S. Iso, D. Karabali and B. Sakita, \PLB{296}{1992}{143}. }
\refdef[TwoDGra]{M. Fukuma, H. Kawai, R. Nakayama, \CMP{143}{1991}{371};\nl
             H. Itoyama and Y. Matsuo, \PLB{262}{1991}{233}.}
\refdef[Fluid]{S. Nojiri, M. Kawamura, A Sugamoto, \MPLA{9}{1994}{1159},
and hep-th/9409164}
\refdef[atlas]{J.~M. ~Figueroa-O'Farrill, J. ~Mas and E.~Ramos.
{\it The topography of $W_\infty$ type algebras}.
\PLB{299}{1993}{41-48}. }
\refdef[Wpart]{E.~ Ramos and J.~Roca. \NPB{452}{1995}{705-723}.
({\tt hep-th/9504071}).}
\refdef[KR]{ V.~G.~ Kac and O.~ Radul, \CMP{157}{1993}{429-457}.}
\refdef[Zoup]{see for example: A.~A.~ Kehagias, P.A.A.~ Meesen and G.~
Zoupanos.\PLB{346}{1995}{262-268}, ~ and references therein.}
\refdef[Kup]{B.~A.~ Kupershmidt,~\CMP{99}{1985}{51}.}
\refdef[OeStra]{W.~Oevel and W.~ Strampp,~{\it Constrained KP hierarchy and
bi-hamiltonian structures},~\CMP{157}{1993}{51}.
}
\refdef[Semenov]{M.A. Semenov-Tyan-Shanskii, {\it What is a classical
R-matrix}, \FAaIA{17}{1983}{259}. -{\it Dressing transformations and
Poisson-Lie
groups}, {\sl RIMS, Kyoto. Univ.} {\bf 21} {\sl ~(1985)~1237}.}
\refdef[NotasJose]{J.M. Figueroa-O'Farrill, {\it $W$-(super)algebras and
(super)integrable hierarchies}. {\sl Unpublished notes, Bonn 1992}.}
\refdef[KuWiYu]{B.A. ~Kupershmidt and G.~ Wilson,~ \Invm{62}{1981}{403}.\nl
F. Yu, ~\LMP{29}{1993}{175}.}
\refdef[KPMiura]{
 J.~Mas and E.~Ramos, {\it The constrained $KP$ hierarchy and the generalised
Miura transformation},~  \PLB{351}{1995}{194} ~({\tt q-alg/9501009}). \nl
J.M. Figueroa-O'Farrill and Sonia Stanciu,~{\it  Lie-Poisson groups and the Miura
transformation},~({\tt q-alg/9501027}).}
\refdef[BaKeKi]{I.~Bakas, B.~Khesin and E.~Kiritsis,~  {\it The Logarithm of
the derivative operator and higher spin algebras of $W_\infty$ type}.
\CMP{151}{1993}{233-243}.}
\refdef[KK]{O.~S. Kravchenko and B.~Khesin, {\it A nontrivial extension of the
Lie algebra of pseudodifferential operators on the circle}.
\FAaIA{25}{1991}{83-85}.}
\refdef[FMR]{J.~M.~Figueroa-O'Farrill, J.~Mas, and E.~Ramos. {\it A one
parameter family
of hamiltonian structures for the KP hierarchy, and a continuous deformation of
the
nonlinear $W_{KP}$ algebra.} \CMP{158}{1993}{17}.}
\refdef[Fernando]{F.~ Mart\'\i nez Mor\'as,  private communication.}
\refdef[Li]{L.~C.~Li and S.~Parmentier, \CMP{125}{1989}{545}.}
\refdef[Gieseker]{D. Gieseker,~{\it The Toda hierarchy and the KdV hierarchy},~
{\tt alg-geom/9509006}.}
\refdef[Aratyn]{H. ~Aratyn, L.A. ~Ferreira, J.F.~ Gomes
       and A.H. ~Zimerman, {\it Toda and Volterra lattice equations from
discrete
symmetries of KP hierarchies},~ \PLB{316}{1993}{85-92},~{\tt hep-th/9307147}.}
\refdef[Bonora]{ L.~Bonora, C.P.~Constantinidis and E.~Vinteler,~
{\it Toda lattice realization of integrable hierarchies},~ {\tt
hep-th/9511172}.}
\refdef[Taka]{K. Takasaki and T. Takebe,  {\it Integrable hierarchies and
dispersionless limit},\nl  {\tt hep-th/9405096}. \nl
 K.~Ueno and K.~Takasaki,~ {\it Toda Lattice Hierarchy }.
Advanced Studies in Pure Mathematics, Vol 4. (1984): Group Representations and
Systems
of Differential Equations, {\sl North Holland} (1984).}
\refdef[Calogero]{
A.P. Polychronakos, ~\PRL{69}{1991}{703}; \nl
L. Brink,~ T.H.~ Hansson and M.A. ~Vasiliev,~\PLB{286}{1992}{109}.\nl
 V. Narayanan and M. Sivakumar.~{\it On the W-algebra in the
Calogero-Sutherland model using the Exchange operators},~{\tt hep-th/9510239}.}
\refdef[JoSo]{J.M. Figueroa-O'Farrill and Sonia Stanciu,~{\it 
Lie-Poisson groups and the Miura transformation},~{\tt q-alg/9501027}.}
\refdef[Feher]{L.Feher, J. Harnand and I. Marshall,~\CMP{154}{1993}{181}.}
%

%%%%%%%%%%%%%%%%%%%%%%%%%%%%%%%%%%%%%%%%%%%%%%%%%%%%%%%%%%%%%%%%%
%
\def\pubblock{ \line{\hfil\twelverm US--FT-32/95}
               \line{\hfil\twelvett q-alg/9512025}
               \line{\hfil\twelverm December 1995}}
%               \line{\hfil\twelverm January 1996}}

%%
\mathstyle
\titlepage
%%
%%%%%%%%%%%%%%%%%%%%%%%%%%%%%%%%%%%%%%%%%%%%%%%%%%%%%%%%%%%%%%%%%%
%%
\title{The algebra of $q$-pseudodifferential symbols and the
$q$-$\WKP^{(n)}$ algebra}
\author{Javier Mas \footnote{$^\sharp$}{\tt
e-mail: jamas@gaes.usc.es}}
\andauthor{Marcos Seco\footnote{$^\flat$}{\tt e-mail:
mseco@fpaxp1.usc.es}}
\Santiago
\abstract{
In this paper we continue with the program to
explore the topography of the space of $W$-type algebras. In the present case,
the starting point is the work of Khesin, Lyubashenko and Roger on the algebra
of $q$-deformed pseudodifferential symbols and their associated integrable
hierarchies. The analysis goes on by studying the associated hamiltonian structures
for which compact expressions are found. The fundamental Poisson brackets yield
$q$-deformations of $\WKP$ and related $\W$-type algebras which, in specific cases,
coincide with the ones constructed by Frenkel and Reshetikhin. The construction
underlies a continuous correspondence between the hamiltonian structures of the Toda
lattice and the KP hierarchies.}
\endtitlepage

\section{Introduction}

The  literature concerning the so called
$\W$ algebras increases as the belief that the $w$ stands for ``wild".
In fact they are wild objects in that they still resist  all efforts to achieve
a clear and unified understanding of their physical meaning or at least of
their geometrical origin. On the other side, the fascination about them stems
from the way they
underlie so many {\sl a priori} disconnected physical and mathematical
constructions: 2
dimensional conformal field theory \[Wn], soliton systems
\[Dickey], vertex-operator and Kac-Moody algebras \[DS], classical and quantum
fluids
\[Hall]\[Fluid], 2-D quantum gravity \[TwoDGra], generalized particle systems
\[RaRo], and a long etcetera.

On the way to taming the $\W$ algebras different proposals have been pursued.
On one hand, in the last
years some effort has been posed in setting up a classification program. It
has been
realized that a natural arena to handle this program is the phase space of
integrable soliton-systems, where very many of the known $\W$ algebras arise
either as
Poisson bracket algebras, or as symmetries of the evolution equations.
On the other hand, searching for an interpretation of $\W$
algebras in physical terms, some simplifications have been produced, yielding
somewhat simpler objects which still preserve many of the distinguishing
features of $\W$
algebras.  Among them, the presence of the Virasoro
subalgebra plays a central role.
Thus for example the ``dispersionless" or ``classical" limit in which
the operator $\d$ is smoothly replaced by a commuting symbol $\xi$ \[Class]
has shed some light  about the geometry of classical
$\W$-morphisms in relation to ``area preserving diffeomorphisms" \[winfty]
and hamiltonian mechanics \[FigRaSo].
\footnote{$^1$}{See also
\[Wgeom] for other interesting proposal.}

Another interesting simplification should occur if we replaced the derivative
$\d$ by the $q$-derivative $\d_q$. The $q$-derivative is in fact a difference
operator, \ie, let $F$ denote the ring of complex valued polynomials in $z$ and
 $z\inv$
($\comps [ z, z^{- 1} ]$) and $q \in \comps$:
$$
\d_q f(z) \equiv {f(qz)-f(z) \over z(q-1)}
$$
$\d$ is recovered in the $\lim_{q\to 1} \d_q = \d$. Using $\d_q$ instead of
$\d$
provides a sort of short distance cuttoff. For this reason it has been widely
investigated in connection with the problem of regulating quantum field
theories
\[Zoup].
In the last years, a few works have been concerned with the issue of the
$q$-deformed
Virasoro an $\W$ algebras; in \[QuantW]\[FResh] we have listed the references
we are aware of, where structures deserving such name have been constructed.
The generic approach in them exploits heavily the use of the $q$-affine
 algebras,
$q$-vertex operators, and a $q$ deformed version the Miura transformation.
The connection of these algebras with integrable systems remained unclear,
until  the recent work of E. Frenkel \[Frenkel]. In this paper it is  claimed
that the $q$-deformed $W$ algebras  constructed in ref. \[FResh] provide
bi-hamiltonian structures for a particular set of differential$-q$-difference
integrable systems, which naturally deserve the name of $q$-deformed KdV
hierarchies.

Our original motivation was to pursue the line of research
developed in \[KLR]. In this work the
central object of study was the Lie algebra of so called $q$-pseudodifferential
symbols $\d_q$, its extensions and contractions, as well as the associated Lax
systems. Actually, the $q$-deformed n-KdV integrable hierarchies defined there
turn out to be the same as those in \[Frenkel], albeit in a different basis.
With respect to this work, ours is somewhat complementary in that we asked
ourselves, first, what are the most general hierarchies that one could
write in terms of Lax operators involving $q$-pseudodifferential symbols and,
second, what are their hamiltonian structures. To perform the analysis, the
unified framework described in \[OeRa] proved to be instrumental. As an
output, a large class of $q$-deformations of classical $W$ algebras are found,
including those of $\WKP$, $GD_n$, or the centrally extended $W_{1+\infty}$. 
In specific cases we find agreement with the results of ref. \[Frenkel].
We also comment on some obstruction found when trying to define a q-deformation of
$W_n$.

This paper is organized as follows: for completeness, sections 2 and 3 are
devoted to the introductory material. In the former one, some
basic notions about the algebra of $q$-pseudodifferential operators are
included; the
later gives an overview of the r-matrix approach to integrable systems. In
both sections we have followed closely the clear expositions of refs. \[KLR]
and
\[OeRa] respectively.

Section 4 is a straightforward application of the machinery of section 3. The
analysis
is performed in a twisted basis $T$, which we refer to as the ``Toda
lattice" basis. In particular, three tri-hamiltonian hierarchies of non-linear
differential-difference equations are found. The Poisson brackets are
explicitly
computed and agree in special cases with those found in \[Frenkel]. One of the
advantages of the present formalism is the possibility of carrying out a
transparent
treatment of reductions. Some of them are investigated at the end of this
section.

Section 5 is a re-elaboration of the previous findings in the basis $\d_q$
introduced in
section 1, and named $q$-KP basis after its direct relationship with the
standard KP basis.
The non-linear infinite dimensional algebra which we obtain and compute is
connected
 with the $\WKP^{(n)}$ algebra \[FMR] in
the limit $q\to 1$; thereafter we name it,
 the $q$-$\WKP^{(n)}$ algebra. Reductions are treated at the end. Of utmost
importance are the reductions of $q$-KP to $q$-KdV. We comment about the
possibility to
obtain several $q$-deformations of the Virasoro algebra within the present
formalism.

Finally, in section 6 we bring the logarithm of the
$q$-differential symbol $\log \d_q$ into the game. We do this by formally
continuing
the order
$n$ of the Lax operator to real values and taking afterwards a suitable limit
$n\to 0$.
The resulting algebra can be considered as a $q$-deformation of the centerful
$W_{1+\infty}$ algebra.

\section{The algebra of $q$-pseudodifferential operators}
 It will be useful to define the ``shift" $ \tau f(z) = f(qz) ,~
\tau^\beta
f(z) = f(q^\beta z),~ \beta\in \comps .$
So, $\d_q$ is a $q$-derivative in the following sense
$$
\d_q (fg) =\d_q (f) g + \tau (f)\d_q (g). \(quleib)
$$
which can be proven by explicit computation.
The actions of $\tau$ and $\d_q$ are not commutative but rather
$q$-commutative, \ie\
$\d_q(\tau(f)) = q \tau(\d_q f)$

\Def<qpedo>
{ An algebra $\pdo_q$ of $q$-pseudodifferential operators is a
vector space of formal series
$$
\pdo_q =
\lbc  A \left({ x,\d_q }\right) =
\sum_{- \infty  }^{ n} u_i (z)\d_q^i \bigm| u_i \in \F \rbc
\(qupedos)
$$
with respect to $\d_q$. The multiplication law in
$\psi DO_q$ is defined by the following rule: $F$
is a subalgebra of $\psi DO_q$ and there are commutation relations ($u\in \F$):
$$
\veqnalign{
\d_q  ~ u &= \left({\d_q u}\right) + \tau (u)\d_q,
\cr
\d_q^{- 1} ~ u &= \sum_{ k \geq 0 } ( - 1)^k q^{-k(k+1)/2}
\left( \tau^{- k - 1} \left( \d_q^k u  \right) \right)\d_q^{- k -1},
\(pasardes) \cr}
$$
}
Each term of the product of two Laurent series in $\d_q$ is found by applying
these rules finite number of times. The formula \(pasardes) is built so that
$\d_q^{-1} ~ \d_q ~ u =\d_q ~ \d_q^{-1} ~ u = u$.
For $q = 1$ these formulas recover the ``classical'' definition of
multiplication law in the algebra of pseudodifferential operators
$\psi DO$.

The commutation rule for $\d_q^n $ (with any integer $n$) and $u(z)$ join
these formulae in one
$$
\d_q^n ~ u =
\sum_{k \geq 0 }^{ }
 \qcomb[n/k]
\left({ \tau^{n - k}
\left({\d_q^k u }\right)
 }\right)\d_q^{n - k},
\(qleibnitz)
$$
where we use the following notation for $q$-numbers and $q$-binomials.
$$
(n)_q ={q^n - 1 \over q - 1}
$$
$$
\qcomb[m/k] =
{(m)_q (m - 1)_q \cdots (m - k + 1)_q
 \over (1)_q (2)_q \cdots (k)_q}.
$$
The $q$-analog of the Leibnitz rule of multiplication of two
$q$-pseudo\-dif\-ferential operators
$A \left({ x,\d_q }\right), B \left({ x,\d_q }\right)$ can be written
as the following operation on their symbols
$$
A \left({ x,\d_q }\right) \,
B \left({ x,\d_q }\right) =
\sum_{ k \geq 0 }^{ }
{1 \over (k)_q ! }
\left({ { d^k \over d\d_q^k}  A}\right) *
\left({\d_q^k B }\right)
\(qsymbolcomp)
$$
where for any complex value of $\alpha$
$$
{d^k\ov d \d_q^k}(f\,\d_q^\alpha) =
(\alpha)_q(\alpha-1)_q\cdots(\alpha-k+1)_q~ f\, \d_q^{\alpha-k}
$$
and the $*$ multiplication of symbols obeys the
following commutation rule for the generators:
$$
 f*\d_q = f\d_q, \qquad \d_q*f = \tau(f) \d_q,
\qquad \d_q^{-1}*f = \tau^{-1}(f)\d_q^{-1}. \(multisymbol)
$$
This follows by a
straightforward verification of the formula \(qsymbolcomp) for the product
$\d_q^n ~ u (z) $, which gives the same answer as \(qleibnitz).

Define the Lie algebra $\G_q$ as the set $\pdo_q$ of all
$q$-pseudodifferential symbols equipped with the commutator bracket
$[A, B] = A ~ B - B ~ A$.

With this setup in mind, it is straightforward to construct a $q$-deformed
analog of the
KP hierarchy. The phase space for this dynamical system is the set $\{ L_q
=\d_q
+ u_1(z) + u_2(z)\d_q^{-1} + u_3(z)\d_q^{-2} + ...\} $, and the equations of
motion adopt
the familiar Lax form
$$
\der{L_q}{t_m} = \comm{L_q}{\left( L_q^m \right)_+} =
\comm{\left( L^m_q \right)_-}{L_q} \(qkplax)
$$
Notice that unlike in the differential case, the potential $u_1(z)$ has a
nontrivial
evolution. This is due to the fact that now the highest degree of the
commutator of two
$q$-pseudodif- ferential operators is the sum of
their respective highest degrees, this being a consequence of the non
conmutativity of
the multiplication of symbols as shown in \(multisymbol).

\section{R-matrix approach to integrable systems}
We recall here the rudiments of r-matrix and refer the interested reader to
the literature \[Semenov]. In this section we shall follow closely the clear
introduction
given in \[OeRa]. A classical
r-matrix on a Lie algebra
$\g$ is a linear map
$\R:\g \to \g$ such that the modified bracket
$$
\comm{a}{b}_\R = \comm{\R(a)}{b}+ \comm{a}{\R(b)}
$$
is a Lie bracket, thus providing a second Lie algebra structure on $\g$.
As was shown in \[Semenov] a sufficient condition for a linear map $\R$ to be
an r-matrix is given by the so-called {\it modified Yang Baxter}  equation
(m-YB($\alpha$) for short).
$$
\comm{\R(a)}{\R(b)}- \R(\comm{a}{b}_\R) = - \alpha \comm{a}{b} \(modiyb)
$$
where $\alpha $ is any real number.
Now let us assume, that in $\g$ there is an ad-invariant
(under the natural Lie bracket $\comm{}{}$ in $\g$) inner product
$\pair{~}{~}:\g\times\g
\to \comps$ under which
$\g$ can be identified with its dual $\g^*$. Immediately we know of a natural
Poisson
structure that lives on
$\comps^\infty (\g^*)$, namely the Lie-Poisson bracket arising from the
modified Lie
bracket $\comm{~}{~}_\R$:
$$
\pbr{f_1}{f_2}_1(L) \equiv \pair{L}{\comm{\R df_1}{df_2}+ \comm{df_1}{\R df_2}
}.
\(liepoisr)
$$
evaluated at a point $L\in \g=\g^*$. This Poisson bracket, termed {\sl linear}
after
its dependence on $L$, is the first of a series of other ``potential" Poisson
brackets.
$$
\veqnalign{
\pbr{f_1}{f_2}_2 &\equiv \pair{L}{~\comm{ \R(L~df_1+df_1~L)
}{df_2}+\comm{df_1}{
\R(L~df_2+df_2~L) }  }\(quadratic)\cr
\pbr{f_1}{f_2}_3 &\equiv
\pair{L}{~\comm{\R(L~df_1~L)}{df_2}+\comm{df_1}{\R(L~df_2~L )} }
\(liepoisrr) \cr }
$$
Using  ad-invariance of the inner product, and the
definition of the adjoint r-matrix as $\pair{\R(a)}{b}= \pair{a}{\R^*(b)}$,
we may encode the above ``potential" Poisson brackets in terms of the
associated
Poisson map $\J$, defined by
$$
\pbr{f_1}{f_2}_s(L) = \pair{J^{(s)}_L(df_1)}{df_2}~~~,~~~~~~ s=1,2,3.
\(hamimap)
$$
as follows
$$
\veqnalign{
      J^{(1)}_L\,(df) &= \comm{L}{\R(df)} + \R^*(\comm{L}{df})  \cr
      J^{(2)}_L\,(df) &= \comm{L}{\R(L~df+df~L)}
+L~\R^*(\comm{L}{df})+\R^*(\comm{L}{df})L\cr
      J^{(3)}_L\,(df) &= \comm{L}{\R(L~df~L)} +L~\R^*(\comm{L}{df})L
\(estructuras)
}
$$
Now the crucial question: for what $\R$ will the above maps define hamiltonian
maps? The findings of \[Li]\[OeRa]  specify that:
\item{$a)$} $J^{(1)}_L$ is hamiltonian for any r-matrix  $\R$ on $g$.
\item{$b)$} $J^{(2)}_L$ is hamiltonian if  $\R$ and its skew-symmetric
combination $\med(\R-\R^*)$ {\it both} satify the m-YB($\alpha$) equation
\(modiyb).
\item{$c)$} $J^{(3)}_L$ is hamiltonian if $\R$ an r-matrix which satisfies
m-YB($\alpha$) equation.

The three maps are related with one-another by simple deformations
$$
\veqnalign{
J^{(2)}_{L+\epsilon 1} &= J^{(2)}_{L} + 2\epsilon J^{(1)}_{L} \cr
J^{(3)}_{L+\epsilon 1} &= J^{(3)}_{L} + \epsilon J^{(2)}_{L} +\epsilon^2
J^{(1)}_{L}
\cr }
$$
where $1$ is the generator of the center in $g$. This, by the way, shows the
compatibility of the  three ``would be" Poisson structures.

The construction of integrable systems, that are hamiltonian with respect to
the above
brackets refers to the existence of a (possibly maximal) set of conserved
functions in
involution. Here, an important piece in the game is played by the set of Casimir
(invariant)
functions, \ie those functions $C\in  \comps^\infty(g^*)$ satisfying
$ad_L^*(C(L)) =
0$ or, equivalently,
$$
  ad_L(dC(L)) = \comm{L}{dC(L)} = 0
$$
If one has a chance to characterize the Casimir functions (in short, the
centralizer of $L \in g$), then a short look at the form of $J^{(s)}$ in
\(estructuras)
reveals that
\item{$(i)$} the associated hamiltonian flows adopt the Lax form
$$
\veqnalign{
\der{L}{t} &= J^{(1)}_L \, (dC) = \comm{L}{\R(dC)} \cr
\der{L}{t} &= J^{(2)}_L \, (dC) = \comm{L}{\R(2L ~dC)} \cr
\der{L}{t} &= J^{(3)}_L \, (dC) = \comm{L}{\R(L^2 ~dC)} \cr
}
$$
\item{$(ii)$} the Casimir functions are in involution.
 For example when $s=1$
$$
\pbr{C_1}{C_2}_1 = \pair{\comm{L}{\R(dC_1)}}{dC_2} =
-\pair{\comm{L}{dC_2}}{\R(dC_1)} = 0.
$$
A particular (partial) solution is given by the
traces of powers of $L$.
$$
C_p(L) \equiv {1\ov k} \Tr (L^p)~,~~~~ dC_p(L) = L^{(k-1)},~~~p=1,2,...
$$
for this particular set of functions, the Lax equations are tri-hamiltonian
$$
\der{L}{t_p} \equiv \comm{L}{\R(L^p)} = J^{(1)}_L(dC_{p+1}) =
J^{(2)}_L (dC_p) = J^{(3)}_L(dC_{p-1})
$$
In some cases ($n$-KdV), $p$ may be a fraction of the order of $L$.

The classification of solutions to \(modiyb) has been achieved partially. A
class of them fall into the following characterization: if $\g= \g_+\oplus
\g_-$ is
a decomposition into Lie subalgebras, denoting by $\P_+$ ($\P_-$) the
projection of
$\g_+$ (resp. $\g_-$) along $\g_-$ (resp. $\g_+$), then $\R= \med(\P_+ - \P_-)$
satisfies
the modified Yang Baxter equation \(modiyb) with $\alpha = 1/4$ since
$\comm{a}{b}_\R$ is easily
calculated to be $\comm{a_+}{b_+} - \comm{a_-}{b_-}$ in obvious notation.

We may give the particular form of \(estructuras) whenever adapted to the
present
situation.\footnote{$^2$}{Hereafter we shall obviate the dependence of
$J^{(s)}_L$ on
$L$, and write simply $J^{(s)}$.}
$$
\veqnalign{
J^{(1)}\, (df) &= ~~~\comm{L}{\P_+ ~df} - \P^*_-\comm{L}{df} \cr
 &= - \comm{L}{\P_- ~df} + \P^*_+\comm{L}{df} \cr
J^{(2)}\, (df) &= ~~~\comm{L}{P_+(Ldf+df~L)}- L~(P_-^*\comm{L}{df}) -
(P_-^*\comm{L}{df})~L \cr
 &= - \comm{L}{P_-(Ldf+df~L)} + L~(P_+^*\comm{L}{df}) + (P_+^*\comm{L}{df})~L
\cr
J^{(3)}\, (df) &= ~~~\comm{L}{P_+(Ldf~L)} - L(P_-^*\comm{L}{df})L \cr
&= - \comm{L}{P_-(Ldf~L)} + L(P_+^*\comm{L}{df})L
\(corchetes) \cr}
$$

Moreover, if $\g_+$ and $\g_-$ are isotropic, then clearly $\R$ is skew-adjoint
with
respect to the inner product $\pair{\R(a)}{b} = \pair{a}{\R^*(b)}$, \ie\  $\R^*
= - \R$.
In this case $\P^*_\pm = \P_{\mp}$
and the 3 structures in \(estructuras) reduce to the following
form $(X\equiv df)$:
$$
\veqnalign{
      J^{(1)}_L(z) &= ~\comm{L}{X_+}_- - \comm{L}{X_-}_+
\cr
      J^{(2)}_L(z) &= ~L(XL)_+ - (LX)_+L =  - L(XL)_- + (LX)_-L  \cr
      J^{(3)}_L(z) &= ~\comm{L}{(LXL)_+}_- - L\comm{L}{X}_+L
\(tresbrackets)  \cr}
$$

\section{ The "Toda lattice" basis}
Let us return to the $q$-deformation of the KP hierarchy that we showed in the
introduction \(qkplax). Define, for $q\neq 1$
$$
\veqnalign{
T &= z(q-1)\d_q + 1 \cr
T^{-1} &= {1\ov z(q-1)\d_q + 1} = \sum_{i=1}^\infty -{(-q)^i\ov
(q-1)^i}z^{-i}\d_q^{-i} \(defdete) \cr}
$$
Any element of $\pdo_q$ of the form \(qupedos) admits a similar expresion in
this
``twisted" basis $$
A = \sum_{-\infty}^n a_i(z) \d^i_q = \sum_{-\infty}^n t_i(z) T^i
\(twistedbasis)
$$
Hence we will be describing the same algebra $\G_q$ in this basis. The
relevant composition law is the following, which can be proven by elementary
manipulations: For
any $f\in \F$
$$
T f = \tau(f) T \(algtau)
$$
in particular $Tz = qzT$. We will use the notation $(Tf)$ to mean that $T$
acts only on $f$, \ie $(Tf) \equiv TfT\inv =\tau(f)$.

The algebraic approach to integrability relies heavily on the existence of an
ad-invariant symmetric bilinear form.
As a step in this direction, a linear functional $\int: f\in \comps$ is defined
satisfying
$\int \tau(f) = \int f$ for all $f\in \F$. In agreement with this requirement,
 we further specify that $\int z^n = \delta_{n,0}$.
A particular realization of this definition is given by the usual Riemann
integration
over $S^1$ of the Fourier basis functions $z^n=e^{in\theta}$, where the action
of
$\tau$ is seen as a shift of $(-i\log q)$ in $\theta$. Also
$\delta(z) = \sum_{k\in \integ} z^k$.

Now, let $A = \sum_i a_i T^i \in \pdo_q$. We define the residue
$\ResT: \pdo_q \to f$ by
$$
\ResT\left( \sum_i a_i T^i\right) = a_0
$$
and the trace $\Tr:\pdo_q \to \comps$ by
$$
\Tr A = \int \ResT A
$$

\Lem<symbili>{ The bilinear form $\pair{~}{~}:\G_q\times \G_q \to\comps$,
given by
$$
\pair{A}{B} = \Tr AB = \int \ResT AB \(pairing)
$$
defines an ad-invariant bilinear symmetric inner product in $\G_q$}

\Pf
By direct computation and use of the defining ``shift" invariance of $\int$ we
find
$$
\veqnalign{
      \Tr AB &= \int \ResT a_iT^i\, b_jT^j  \cr
            &= \int  a_i\tau^i(b_{-i})
            ~= \int  b_{-i}\tau^{-i}(a_i) \cr
            &= \int \ResT b_j\tau^{j}(a_i)~ T^{j+i}
            ~= \int \ResT b_jT^j \, a_iT^i \cr
          &= \Tr BA
{}~~~~~~~~~~~~~~~~~~~~~~~~~~~~~~~~~~~~~
{}~~~~~~~~~~~~~~~~~~~~~~~~~~~~~~~~~~~~~~\QED \cr }
$$
We would like to stress that this bilinear product is the same (up to factors of
$q$) as the one defined in \[KLR], as we shall show in section $5$. The previous
lemma is fully equivalent to theorem $3.3$ in that reference. With respect to this
inner product, the adjoint of $\tau$ is 
$\tau^*=\tau\inv$, \ie\  $(T^*f) = (T\inv f)$.

Let us investigate the possible splittings of the form
$\G_q = \G_1 \oplus \G_2$, where $\G_1$ and $\G_2$ are Lie subalgebras. In
view of the
generic (graded) commutation relations
$$
\comm{t_i T^i}{t_j T^j} = (t_i\tau^i(t_j)- t_j\tau^j(t_i)) T^{i+j}
$$
  we find only three possibilities as follows

\item{1.} ($\underline{\sigma=-1}$),  $~~~ \G_q
= \G_{\geq 0}\oplus \G_{\leq -1}$
$$
\G_{\geq 0}  \equiv \{~\sum_{i \geq  0} t_i(z) T^i   ~ \} ~~~~;~~~~~
\G_{\leq -1} \equiv \{~\sum_{i \leq -1} t_i(z) T^i ~\}
$$
\item{2.} ($\underline{\sigma=+1}$), $~~~\G_q=\G_{\geq 1}\oplus \G_{\leq
0}$:
$$
\G_{\geq 1}  \equiv \{~\sum_{j \geq 1} t_j(z) T^j ~\} ~~~~;~~~~~
\G_{\leq 0}  \equiv \{~\sum_{j \leq 0} t_j(z) T^j  ~ \}
$$
\item{3.} ($\underline{\sigma=~~0}$),  $~~~~\G_q = \G_{0^+}\oplus \G_{0^-}$
$$\veqnalign{
\G_{0^+} &\equiv \{\sum_{k\geq 0} t_k(z) T^k, ~t_0\in z\comps[z]    \} ~~;~~
\cr
\G_{0^-} &\equiv \{\sum_{k\leq 0} t_k(z) T^k,~ t_0\in z\inv\comps[z\inv] \} \cr
}
$$
\Rmk<manin>{ As mentioned in \[KLR],
the interest of the last case comes from the fact that, relative to the inner
product defined
in \(pairing), it is the only one where $\G_{0^\pm}$ are isotropic.
 Hence $(\G,\G_{0^+},\G_{0^-})$ is a Manin triple and $(\G,\G_{0^-})$ a Lie
double.
}
\subsection{The fundamental Poisson brackets}

In order to define the phase space where our dynamics will take place,
let
$$
L = \sum_{i=m}^{n} t_i(z) \,T^i \(opelaxt)
$$
where $n,m\in \integ$ and $n>m$,
We regard any of these difference  operators as ``points" on a manifold
$\M^{(n,m)}_T$.
The dynamics is governed by differential$-q$-difference equations derived from the
usual Lax system
$$
\der{L}{t_p} = \comm{L}{\left( L^p \right)_+} =
\comm{\left( L^p \right)_-}{L} \(laxecu)
$$
which is manifestly consistent for $L$ of the form given in \(opelaxt)
Here $\pm$ refers to any one of the $\sigma=0,\pm1$ splittings defined above.

In the case of $\sigma=-1$ with $(n,m)=(1,-\infty)$
this system is none other than
the simplest version of the Toda lattice hierarchy involving one set of time
parameters
\[Taka]\[Bonora]. Indeed, in these works the Toda lattice hierarchy is
formulated in terms of a Lax operator of the form 
$$
L = e^{\d} + \sum_{n=0}^\infty u_{n+1} e^{-n\d}
$$
which involves the difference operator  $e^{\d}$ acting as
$e^{n\d} u_i(x)= u_i(x+n)e^{n\d}$. The isomorphism between both
formulations is made patent after identifying $u_i(x)$ with $t_i(z=q^x \zeta)$
where $\zeta\in\comps$ is any fixed complex number.

Remark that only for $m=0$ and the splitting
$\sigma=-1$, or $n=0$ and $\sigma=+1$ equations \(laxecu) are empty since in this
case the commutator vanishes identically. The non-trivial flows may come then from
fractional powers of $L$ \footnote{$^3$}{This issue is more delicate than in
the usual context of KdV, and requires a careful definition of the ring of
functions \[Frenkel]\[KLR]. We will not dwell here with this aspect, certainly
important
from the point of view of  integrability.}. For example, let
$n=N$ and $m=0$, then
$$
\der{L}{t_p} = \comm{L}{\left( L^{p/N} \right)_+} =
\comm{\left( L^{p/N} \right)_-}{L}
$$
are non-trivial differential difference equations as long as $p$ is not a
multiple of
$N$.
An analogous way to characterize these flows is to consider an operator of the form
$L\in \M_T^{1,-\infty}$ of the form  $ L = t_0 T + t_1 + t_2 T^{-1}+....$,
constrained to satisfy  $L^{N}_- = 0$.

On $\M^{(n,m)}_T$ the (linear) functionals of interest have the form
$f_X(L) = \Tr LX$
with
$$
X= \sum_{j=m}^n T^{-j} x_{j}(z) .
$$
 Clearly $f_X$ adopts
the form of an euclidean scalar product $f_X = \int \sum_{i=m}^n t_i x_i$ .
Defining the gradient $d: \F \to \G_q$ by
$$
\pair{df}{\delta L} \equiv \left.
\der{}{\epsilon} f_X (L+\epsilon\delta L) \right\vert_{\epsilon=0}
$$
it turns out that $df_X(L) = X$.

We are interested in the fundamental Poisson brackets among the fields
$t_i(z)$.
Since  the Poisson maps $J^{(s)}(df_X)$ are linear in $df_X = X$ we may expand
$$
J^{(s)}(X) = \sum_{i=m}^n\sum_{j=m}^n (J^{(s)}_{ji} \, x_i)(z) \,T^j
\(linmap)
$$
where $J^{(s)}_{ij}$ is some function of q-difference operators $T$.
Plugging this back in \(hamimap) we obtain
$$
\pbr{f_X}{f_Y} = \int (J_{ij} \, x_j) y_i = \int x_i(J^*_{ji} \, y_j)
= -\int x_i(J_{ij}\, y_j)
$$
where the last equation follows from the antisymmetry $\pbr{f_X}{f_Y}=
-\pbr{f_Y}{f_X}$ which implies that $J^*_{ji} = - J_{ij}$. Finally, comparing
this
expression with
$$
\pbr{f_X}{f_Y}_r = \int\, x_i(z)\int\, \pbr{u_i(z)}{u_j(w)}_s\, y_j(w)
$$
shows that
$$
\pbr{u_i(z)}{u_j(w)}_s = -(J^{(s)}_{ij}(z)\, \delta(z/w)) \()
$$
where $\delta(z/w) = \sum_{j\in\integ} (z/w)^j$, and the operators $J_{ij}$ 
act at $z$.

It is time to analyse in detail the potential Poisson structures on
$\M^{(n,m)}_T$.
We will do this by taking into consideration, case by case,
the three possible splittings of
$\G_q$: $\sigma= 0,\pm1$. Notice that for all cases, the linear and cubic
brackets
in \(liepoisr) and \(liepoisrr) define Poisson brackets,
since $\R = \med(\P_+ - \P_-)$, with $\P_\pm$ in each case the relevant
projection
operators, yields automatically an r-matrix obeying the m-YB(${1\ov 4}$)
equation \(modiyb). Therefore, further analysis is only required for the quadratic
bracket.
\item{$1.$} ($\underline{\sigma=-1}$): $~~ \G = \G_{\geq 0}\oplus \G_{\leq -1}$

This splitting is, among the three, the most analogous to the one of
the standard KP hierarchy. Notice however the important difference: now the
subalgebra
$\G_{\geq 0}$ is {\it  not isotropic}, and in consequence the r-matrix
is not anti-selfadjoint.Thus, whether the ``antisymmetric"
combination
$\med(
\R-\R^* )$ satisfies the  m-YB(${1\ov 4}$) equation as well, must be checked
independently. In more concrete terms, let
$$
\R = \med(\P_{\geq 0} - \P_{\leq -1}).
$$
In view of the definition of the inner product \(pairing)
$\R^* = \med (\P_{\leq 0} - \P_{\geq 1})$, and therefore
$$
\med( \R-\R^* ) = \med(\P_{\geq 1} - \P_{\leq -1}) . \(antisymerre)
$$
It follows from an easy calculation that this  linear map satisfies \(modiyb)
with $\alpha=1/4$ as
well. Hence all three brackets in \(liepoisr) and \(liepoisrr) are Poisson
brackets.
Using the general formula \(corchetes) we may particularize to the present case
and get $(X\equiv df)$
$$
\veqnalign{
J^{(1)} ~(X) &= ~~~\comm{L}{X_{\geq 0}}_{\leq 0} - \comm{L}{X_{\leq -1}}_{\geq
1}
\(estruno)\cr
J^{(2)} ~(X) &=~~\,2 L(XL)_{\geq 0} - 2(LX)_{\geq 0} L + L\,\ResT(\comm{L}{X})+
\ResT(\comm{L}{X})\, L  \cr
&= - 2L(X~L)_{\leq -1} + 2(LX)_{\leq -1} L + L\,\ResT(\comm{L}{X})+
\ResT(\comm{L}{X}) L
\(estrudos)\cr
J^{(3)} ~(X) &= ~~~\comm{L}{(LXL)_{\geq 0}} - L\,\comm{L}{X}_{\geq
1}\,L
\(estrutres) \cr}
$$
  A word about consistency. Concerning the {\it first structure}
\(estruno) \footnote{$^4$}{For the particular case of $\M^{(n,0)}_T$ 
this expression also appears in \[Frenkel].} in order that
$J^{(1)}(X)$ describes a deformation of $L$ we must demand that
$n\geq 0 $ and $1\geq m$. This may be seen from \(estruno)
written in the following form
$$
J^{(1)} ~(X) = \comm{L_{\leq 0}}{X_{\geq 0}}_{\leq 0} - \comm{L_{\geq
2}}{X_{\leq
-1}}_{\geq 1}.
$$
The first term on the right hand side bounds the lowest order of 
$J^{(1)}(X)$ to be higher or equal to the lowest order of $L$. The second
one bounds the highest order of the map to be strictly lower than the highest
order of $L$. The constraints in $n$ and $m$ come from the projectors
outside the commutators.

It is remarkable that for the {\it second structure} $J^{(2)}$ there is not such a
restriction and $J(X)$ parameterizes a deformation of $L$ for any $(n,m)$. 

The case of the cubic or {\it third structure} is more extreme. From \(estrutres)
is evident that
$\M^{(n,m)}_T$ is not an invariant subspace under the action of
$J^{(3)}$, unless
$(n,m)=(\infty,0),~ (0,-\infty)$ or $(\infty,-\infty)$. For this reason we do
not
consider this algebra to be of much interest and we will not write down its
explicit
form. Hereafter we shall restrict our attention to the two other structures.

In terms of the fundamental Poisson brackets we obtain for $J^{(1)}$ the
following
difference operators $J^{(1)}_{ij}$
$$
J^{(1)}_{ij} = t_{i+j} T^i - T^{-j} t_{i+j}
 \(estrunocompo)
$$
if either $i,j\geq 1$ with $n\geq i+j$ or, up to a sign, when $0\geq i,j$
with $i+j\geq m$. In all other cases $J^{(1)}_{ij} =0$ .
For the sake of comparison with similar results in the
literature \[FResh]\[QuantW]\[Frenkel], we may write down the Poisson
brackets explicitely.
$$
\pbr{t_i(z)}{t_j(w)}_1 = -t_{i+j}(z)~ \delta\left({q^iz\ov w} \right)
                       + t_{i+j}(w) ~\delta\left({z\ov q^jw} \right)
\(estrunocompa)
$$
with the same set of restrictions upon the indices $i,j$. This expression
exhibits the
splitting of the linear Poisson bracket algebra in 2 graded subalgebras spanned
by either
positive or non-positive values of $i,j$.

 For $J^{(2)}_{~ij}$ an
analogous computation yields
$$
J^{(2)}_{~ij} = 2\sum_{k=\max (m,i+j-n)}^{\min (n,i)} \left(
t_k T^{k-j} t_{i+j-k} - t_{i+j-k} T^{i-k} t_k \right) +
t_i~(1+T^i)(1-T^{-j})~t_j
$$
or, again
$$
\veqnalign{
\pbr{t_i (z)}{t_j (w)}_{(2)}& =~~~ 2
\sum_{k=\max (m,i+j-n)}^{\min (n,i)}  \left(
 t_{i+j-k}(z) t_k(w)~ \delta \left({q^{i-k} z\ov w}\right)\right.
\cr
& \left. -t_k(z)  t_{i+j-k}(w) ~\delta\left( {z\ov q^{j-k}w} \right) \right)
  - t_i(z)t_j(w) \sum_{l\in \integ} \left({z\ov w}\right)^l(1+q^i)(1-q^{-j})
\cr
}
$$
\item{$2.$} ($\underline{\sigma=+1}$): $~~~ \G_q = \G_{\geq 1}\oplus \G_{\leq
0}$

This situation is symmetric with respect to the one above. Notice that at the
level of the algebra, this splitting transforms into the previous one upon the
substitution
$T\to T\inv$. Therefore, the formulas obtained from \(estruno)-\(estrutres) 
can be
adapted to the present case by a simple replacement $q\to q\inv$ and
$m\leftrightarrow
n$.
\item{$3.$} ($\underline{\sigma=~~0}$): $~~~\G_q = \G_{0^+}\oplus \G_{0^-}$

This is the standard case of a Lie bialgebra. The three maps in
 \(tresbrackets) automatically define Poisson brackets.
The fundamental ones are a slight modification of the ones above, and involve
an
additional operator $\p_\pm$, which projects any element $f\in F$ into its
Taylor
and Laurent parts respectively, \ie\ $\p_+ z^m= z^m$ iff $m\geq 1$ and zero
otherwise,
and viceversa. $\p_+$ and $\p_-$ are mutually adjoint with respect to the
inner product defined with $\int$ and commute with $T$.

As before, the linear structure is a direct sum of two subalgebras, spanned by
the fields $t_0^+\equiv \p_+t_0$ and $\{t_i,~i=1,...,n\}$ on ones side, and
$t_0^-\equiv \p_-t_0$ and $\{t_i,~i=-1,-2..,m\}$ on the other.
Therefore as long as $i,j\geq 1$ but $n\geq i+j$
$$
J^{(1)}_{ij} = -(t_{i+j}T^i- T^{-j} t_{i+j}) \()
$$
The same expression with opposite sign holds if $-1\geq i,j$ with
$i+j\geq m$. Finally
$$
J^{(1)}_{0j} = -\Theta(j-1)\, \p_+(T^{-j}-1)t_j
{}~+~\Theta(-j-1)\,\p_-(T^{-j}-1)t_i
\(abracero)
$$
and $J^{(1)}_{i0} = - J^{(1)*}_{0i}$. In all other cases $J_{ij}=0$. In
formula
\(abracero)
$\Theta$ stands for the usual step function $\Theta(x)=1$ iff $x\geq 0$ and $0$
otherwise. The quadratic brackets are computed along the same lines:
$$
J^{(2)}_{~ij} = 2\sum_{k=\max (m,i+j-n)}^{\min (n,i)} \left(
t_k T^{k-j} t_{i+j-k} - t_{i+j-k} T^{i-k} t_k \right) + 2 t_i(1-T^{i-j})\p_-
t_j
$$
\hyphenation{pseu-do-dif-fe-ren-tial}
\subsection{Some reductions}

Let us  focus on the $\sigma=-1$ splitting (the case $\sigma=+1$ follows a
symmetric pattern).
 Remark that as, far as the Lax equations are concerned, the field $t_n$ is not
dynamical. Let $\wt\M^{(n,m)}_{T}$  represent the submanifold of
$\M^{(n,m)}_T$
defined by the constraint $t_n=1$ (or any constant).
From $J^{(1)}$ in \(estruno)
we observe that, as long as $n\geq 1$, the highest positive order of $J^{(1)}$
is $n-1$
and therefore the hamiltonian map is automatically tangent to the constraint
submanifold.  When $m=0$ this is also the case for a similar contraint
on the lowest field $t_0=$ constant;
indeed \(estruno) shows that in this case the contribution of $J^{(1)}(X)$
to order zero is $\comm{L_0}{X_0} = 0$. In few words, both constraints are
first class, and the Poisson brackets are defined by simple restriction of
\(estrunocompa).

For $J^{(2)}$, things are more involved. Notice in fact from the expression
\(estrudos), that the highest order of $J^{(2)}(X)$ is $n$, \ie\ the same as
that of $L$. Therefore, in order to define Poisson brackets on
$\wt\M^{(n,m)}_{T}$ we would follow the standard prescription for second
class constraints due to Dirac. However, instead of plugging here the formula
of the Dirac brackets we will pause briefly to describe how they appear in our
context. Given  the projection map, such as
$\M^{(n,m)}_T \to \wt\M^{(n,m)}_T$ that sets $t_n\to 1$, at each point $L$,
the induced projection of vectors on the tangent subspace is unique. This is
not  so for 1-forms. To see this notice that if we want to compute Poisson
brackets of functions $f,g$ on $\wt\M^{(n,m)}_{T}$ via \(estrudos) we first
need to extend them to $\M^{(n,m)}_{T}$. This extension, being non-unique, 
renders the component $x_n=\delta f/\delta t_n$  in the gradient 
$$
df(L)=\sum_{k=m}^{n} T^{-k}x_k
$$ 
undefined. Therefore some additional structure is required in order to
specify the cotangent subspace. Since we have
$J^{(2)}$ at hand, a map from 1-forms to vectors, we may fix this ambiguity
by demanding that the associated hamiltonian vector fields be tangent to  
$\wt \M^{(n,m)}_T$.
In other words, we fix $x_n$ by the requirement that
$J_2(df(L))$ should have no term of order $n$.
 This form of computing the algebra is fully equivalent to the Dirac bracket
prescription  as we show next \[NotasJose]. The demand that
$J(z)$ should stay tangent to the constraint manifold implies for $df(L)$
that 
$$
\sum_{j=m}^n J_{nj} \, x_j  = 0
$$
and this may be solved for  $x_n = -\sum_{j=m}^{n-1} J^{-1}_{n n}\,
J_{n j}\, x_j$. Plugging this back into \(linmap) we have
$$
J(z) = \sum_{i,j= m}^{n-1} (\wt J_{ij}\, x_j) T^i
$$
where
$$
\wt J_{ij} = J_{ij} - J_{in}\,  J^{-1}_{nn} \, J_{nj},~~~~~~
i,j=1,...,n-1.
\(dirbra)
$$
are the corresponding Dirac brackets on the constraint surface. For the
explicitly
reduced brackets we find a non-local expression as follows
$$
\veqnalign{
\wt J^{(2)}_{ij} = 2\sum_{k=\max (m,i+j-n)}^{\min (n,i)} &\left(
t_k T^{k-j} t_{i+j-k} - t_{i+j-k} T^{i-k} t_k \right)  + 2~ t_i~
{(1-T^{i-n})(1-T^{-j})\ov (1-T^{-n})}~t_j \cr}
\(brackeredudos)
$$
Indeed, the interest in the reduction $t_n=1$
stemmed from the fact that the Lax flows \(laxecu) stabilize this constraint.
Likewise, if $m=0$ the Lax equation for $t_0$ is trivial, hence we may 
want to set it also to a constant. However, in contrast to the previous
case, the contribution of
$J^{(2)}(X)$ to order zero is $2L_0(XL)_0 -2(LX)_0L_0+2L_0\comm{L}{X}_0=0$; in
other words, for all $j$, $J^{(2)}_{0j}$ vanishes and therefore this constraint
is first-class and
does not lead to any modification of the algebra.
If we put $m=0$, \(brackeredudos) is equivalent to formula (3.6) in 
\[Frenkel]. 

\section{The $q$-KP basis}
We recall that our main purpose is to construct a $q$-deformation of the
algebra $\WKP$.
For this reason it will be interesting to reformulate the findings of the
previous
section in terms of the basis $\d_q$, \ie
$$
\d_q = {1\ov z(q-1)}\,(T-1) \(cambioq)
$$
Written in this basis, the limit $q\to 1$  should yield directly
 $\WKP^{(n)}$ in \[FMR].
  We recall here the relevant formulae for the change of basis.
$$
\veqnalign{
T &= z(q-1)\d_q + 1 \cr
T^{-1} &= {1\ov z(q-1)\d_q + 1} \equiv - \sum_{i=1}^\infty {(-q)^i\ov (q-1)^i}
z^{-i}\d_q^{-i}
\(defdeteq) \cr}
$$
 These imply in particular, that the
phase space $\M_T^{(n,m)}$ will be coordinatized now by $q$-pseudodifferential
operators  $L$, of the form
$$
L = \sum_{j=-\infty}^n u_j(z) \d_q^j .\(qulaxq)
$$
$(m+n)$ being still the number of degrees of freedom. Yet the manifold of all
$q$-pseudodif-ferential operators of the form \(qulaxq), which we will denote
by
$\M_{\d_q}^{(n)}$, is much bigger than
$\M_T^{(n,m)}$. Rather we have that the set of  all these spaces
$\{\M_T^{(n,m)},~m=1,2,3...\}$ is dense in $\M_{\d_q}^{(n)}$.

Notice that \(defdeteq) involves a specific choice of the expansion point,
namely around $\d_q =\infty$. Other choices may lead to different $W$-algebras.  
In the
$q$-KP basis we may grade $\G_q$ by the scaling dimension: if $z$ has degree 
$-1$, $\d_q$ will have $+1$ and we may make $L$ homogeneous of a certain degree,
$n$,  by further assignment of degree $j$ to $u_j$. This gives us a chance to look
for a
$q$-deformation of the Virasoro algebra in the subalgebra spanned by the
counterpart of the energy momentum tensor  (the field $u_{2}$ in the context of the
classical
$\WKP$ algebra), which will be a particular $x$-dependent combination of
various fields in the Toda basis (where the grading was a different one).

In order to study the hamiltonian structures we have to re-define the residue
and trace functionals in the new basis.  The point is the following; let
$$
L(T)_{\geq 0} =\sum_{i=0}^n t_iT^i= \sum_{i=0}^n u_i \d_q^i =
L(\d_q)_{\geq 0}
$$
where in each case the
projection is performed with respect to the relevant basis.
Making use of \(cambioq) we may write $t_0$ in terms of $u_i$:
$$
t_0(u_i) = (-1)^m{u_m\ov z^m(q-1)^m} + (-1)^{m-1}{u_{m-1}\ov
z^{m-1}(q-1)^{m-1}} +...-
{u_1\ov z(q-1)} + u_0 \(uveu)
$$
If $L(T) \equiv t_iT^i$, $t_0=\ResT L(T)$. How can we manage extract $t_0(u_i)$
out of
$L(\d_q)$ as given by the right hand side of \(uveu)?
Notice that we may take advantage of the fact that the projections (in the
respective
basis) $(~)_{\geq 0}$ and $(~)_{\leq -1}$ commute with the change of basis
$T\leftrightarrow \d_q$, hence
$$
(L(T))_0 = (L(T)T\inv)_{-1} = (L(T)_{\geq 0}T\inv)_{-1} =
 {z(q-1)\ov q} (L(\d_q)_{\geq 0}T\inv(\d_q))_{-1}
$$
In the last expression
$T\inv(\d_q)$ stands for the second relation in  \(defdeteq). Thus
$$
t_0(u_i) = {z(q-1)\ov q} \,
\ResD (L(\d_q)T\inv).
$$
where we have introduced the symbol $\ResD a_i\d_q^i \equiv
a_{-1}$.
Concerning the ad-invariant symmetric bilinear product, we find the same
expression
that was considered in ref \[KLR] (modulo a constant factor)
$$
\pair{A}{B}  =\pair{B}{A}  = \int \ResT\, AB =
 {q-1\ov q} \,\int z \, \ResD ABT\inv
 \(scapro)
$$
In this basis the natural integral functional is $\int_{-1}\equiv \int z$
which, in
spite of not being scale invariant $\int_{-1} \tau(f) = q\inv \int_{-1} f$
satisfies the
desirable property that $\int_{-1} (\d_qf) = 0$. With respect to the above
inner
product, the adjoints of $\tau$ and $\d_q$ are easy to compute, yielding
$$
\tau^* = {1\ov q} \tau\inv~~~;~~~\d_q^* = -\d_q^*\, \tau\inv .\(adjuntos)
$$
For later use we shall introduce the following compact notation:
$$
\Omega(A) \equiv {z(q-1)\ov q}~ \ResD(A)
$$
Next, we must characterize the three possible splittings of $\G_q$ ($\sigma =
0,\pm 1$) in
the $\d_q$ basis:

 The untwisted basis is naturally adapted to the case $\underline{\sigma= -1}$
$$
\G_{\geq 0}(T) = \G_{\geq 0}(\d_q)~~~;~~~\G_{\leq -1}(T) = \G_{\leq -1}(\d_q).
$$
$\underline{\sigma = +1}$ looks a little bit more contrived
$$
\veqnalign{
\G_{\geq 1}(T) &= 
\tilde\G_{\geq 0}(\d_q) \equiv \{L=\sum_{j=0}^m u_j\d_q^j~|~
\ResD (LT\inv)= 0 ~\} \cr
\G_{\leq 0}(T) & = \G_{\leq 0}(\d_q) = \{L=\sum_{j=-\infty}^0 u_j\d_q^j~
 ~ \} \cr
}
$$
Lastly, the characterization of $\underline{\sigma= 0}$ in the untwisted basis
makes this splitting very unnatural
$$
\veqnalign{
\G_{0^+}(\d_q)  &\equiv \{ L= \sum_{j=0}^m u_j\d_q^j~|~
z\ResD (LT\inv) \in z\comps [z] ~\} \cr
\G_{0^-}(\d_q)  &\equiv \{ L= \sum_{j=-\infty}^0 u_j\d_q^j~|~
u_0 \in z^{-1}\comps [z\inv] ~\} \cr
}
$$
We want to consider again the
Poisson maps \(corchetes). Now in order to compute the analog of
\(estruno)-\(estrutres)
we have to say what the relevant projection operators are. From the form of the
scalar
product \(scapro) it is clear that
$$
\veqnalign{
P_{\geq 0}L = L_{\geq 0 }  ~~~&;~~~ P^*_{\geq 0}L = (LT\inv)_{\leq -1}T \cr
P_{\leq -1}L = L_{\leq -1} ~~~&;~~~ P^*_{\leq -1}L = (LT\inv)_{\geq 0}T
\(duales)\cr
}
$$
We may simplify these expressions, reminding that the projections $(~)_{\geq
0}$ and
$(~)_{\leq -1}$ commute with the change of basis $T\leftrightarrow \d_q$.
So for $L = t_i T^i = u_j \d_q^j$
$$
\veqnalign{
 P^*_{\geq 0}L &= (L(\d_q)T\inv(\d_q))_{\leq -1} T =
(L(T)T\inv)_{\leq -1}T = (L(T))_{\leq 0} \cr &= (L(T))_{\leq -1} + t_0
 \cr & = (L(\d_q))_{\leq
-1} +\Omega(L) \cr }
$$
where we made use of \(defdeteq). Similarly
$$
\veqnalign{
 P^*_{\leq -1}L &= (L(\d_q)T\inv)_{\geq 0}T =
(L(T)T\inv)_{\geq 0}T = L(T)_{\geq 1}T\inv T \cr & = L(T)_{\geq 0} - t_0
\cr  &= L(\d_q)_{\geq 0} - \Omega(L) \cr }
$$

With these results the antisymmetric part of the r-matrix is
$$
\med (\R-\R^*) = \med (P_{\geq 0} - P_{\leq -1} - \Omega) = \R - \med \Omega
\(antisymerreq)
$$
It is not evident that this expression also satisfies
the m-YB(${1\ov 4}$) equation. However an explicit
computation shows that the
only non-vanishing contribution to \(modiyb) has the form $\Omega(\comm{a_{\geq
0}}{b_{\geq 0}})$, which vanishes. An easy way to convince oneself of this fact
is that
when written in the
$T$ basis this is  $\ResT\comm{a(T)_{\geq 0}}{b(T)_{\geq 0}} = 0$.

With all this information, it is now an easy exercise to find the explicit
expressions
for
\(corchetes)  as adapted to the present case.
% The result was to be expected as it
%coincides with
%\(estruno)-\(estrutres) after replacing the functional $\ResT$ by its
%%counterpart
%in the untwisted basis  $\Omega \equiv z(q-1)/q~ \ResD$.
 Let again $X\equiv df$:
$$
\veqnalign{
J^{(1)}(X) &= ~~~\comm{L}{X_{\geq 0}}_{\leq -1} - \comm{L}{X_{\leq -1}}_{\geq
0}
+\Omega(\comm{L}{X}) \(elprimer)
\cr
J^{(2)}(X) &=~~\,2 L(XL)_{\geq 0} - 2(LX)_{\geq 0} L +
L~\Omega(\comm{L}{X})+
\Omega(\comm{L}{X})~ L  \cr
&= - 2L(XL)_{\leq -1} + 2(LX)_{\leq -1} L + L~\Omega(\comm{L}{X})+
\Omega(\comm{L}{X})~ L  \(elsegun)
\cr
J^{(3)}(X) &= ~~~\comm{L}{(LXL)_{\geq 0}} - L\comm{L}{X}_{\geq 0}L
+L\Omega(\comm{L}{X})L
\(corchetesd) \cr}
$$
Notice that as compared with the analogous expressions for the $\WKP$ algebra
\[FMR], the ones above present additional terms which vanish in the limit $q\to
1$. However these terms are not active whenever
$f$ is a Casimir function, and hence, in particular for the Lax-hamiltonian
flows.

In order to compute the algebra of fundamental Poisson brackets we have to
describe the manifold and the class of funtionals for  which $J^{(i)},~i=1,2,3$
describe tangent maps.
We will work on $\M^{(n)}_{\d_q}$ whose points are parameterized as
$$
L^{(n)}=  \sum_{i=0}^\infty u_i \d^{n-i} ~~~~~(n\in \integ)
\(defilaxi)
$$
Accordingly, in order to define linear functionals of the form  $ f_X =
\int_{-1} u_i x_i$ as \nl $f_X =  \int_{-1} \ResD LXT\inv$ our gradient 1-forms
will be
$q$-pseudodifferential operators of the form
$$
X \equiv df_X = \sum_{j=0}^{\infty} \d_q^{j-n-1} ~ T~  x_j.
$$
After a straightforward computation, we list the full set of fundamental
brackets for $J^{(1)}$ as  follows: first we have that  for all $j$:
$
J^{(1)}_{00}=J^{(1)}_{0j}=J^{(1)}_{j0}=0
$
 If $\underline{i,j\geq n+1}$:
$$
\veqnalign{
J_{ij}^{(1)} =&
\sum_{k=0}^{i+j-n-1}\qcomb[i-n-1/k] q^{k(k+1)\ov 2}
\left( q^{n-i-1}(q-1)~u_{i+j-n-k}~(-\d_q)^k ~x + {1\ov q}
u_{i+j-n-k-1}~(-\d_q)^k~
\right) T^{n-i}
\cr
&~~~
-~\sum_{k=0}^{i+j-n}\qcomb[j-n/k]{(q-1)\ov q}~T^{j-n-k}~\d_q^k~ u_{i+j-n-k}~x
\cr
&~~~~~~~
-~\sum_{k=0}^{i+j-n-1}\qcomb[j-n-1/k]{1\ov q}~T^{j-n-k-1}~\d_q^k~ u_{i+j-n-k-1}.
\(mayorquene)\cr}
$$
If however $\underline{1\leq i,j \leq n-1}$ the same expression \(mayorquene)
is valid with the opposite sign. Finally when $\underline{j=n}$:
$$
\veqnalign{
J_{in}^{(1)} =&
\left\{ \sum_{k=0}^{i-n-1}\left(\qcomb[i-n-1/k]q^{{k(k+1)\ov 2}+n-i-1}(q-1)
{}~u_{i-k}
\,(-\d_q)^k\,x\,T^{n-i}\right) - {1\ov q}(q-1)~ x\,u_i\right\} \Theta(i-(n+1))
\cr
&~~~
+~\left\{ -\sum_{k=0}^{i-1} \qcomb[i-n-1/k]q^{k(k+1)\ov 2} \left(
q^{n-i-1}(q-1)~ u_{i-k}(-\d_q)^k~x + {1\ov q}~ u_{i-k-1}(-\d_q)^k\right)\right.
\cr
& ~~~~~~~~
\left.+ {(q-1)\ov q}\,x\,u_i +\sum_{k=0}^{i-1}\qcomb[-1/k]{1\ov q}~T^{-k-1}
\d^k_q \, u_{i-k-1}  \right\} \Theta(n-i).
\cr
&~~~~~~~~~~~~~~~~
+~{x(q-1)\ov
q}\left(~u_n-\sum_{k=0}^n\qcomb[-1/k]~T^{-k}~\d_q^k~u_{n-k}~\right)
 \delta_{i,n} \(jigene) \cr
}
$$
  The rest of the brackets can be
computed making use of the identity $J_{ij}=-J_{ji}^*$, and \(adjuntos)

Concerning reductions, as long as
$n\geq 1$, the highest order of $J^{(1)}(X)$ is $n-1$ and thereafter its action
is
tangent to the submanifold defined by $u_0=$  constant.
Again this reduction is therefore first-class.
Two other consistent reductions of $L$ are of the form $L= L_{\geq 0}$ with
$
J^{(1)}(X) = -\comm{L_{\geq 0}}{X_{\leq -1}}_{\geq 0} +
\Omega( \comm{L_{\geq 0}}{X_{\leq -1}})
$
or $L=L_{\leq 0}$, in which case
$
J^{(1)}(X) = \comm{L_{\leq 0}}{X_{\geq -1}}_{\leq -1}+
\Omega(\comm{L_{\leq 0}}{X_{\geq -1}}) .
$
The relevant explicit form of the Poisson brackets can be obtained in each
case from \(elprimer) after suitably setting to zero the corresponding fields
$u_i$ and its duals $x_i$.

Written in this basis, the formula \(mayorquene) exhibits a nested sequence
of subalgebras $N=1,2,...$, spanned by $\{u_{n+N+k},~k=0,1,2,...\}$. In the
continuum limit $q\to 1$ these contract to the nested set truncations of the
centerless $W_{1+\infty}$ algebra known as $W_{-N+\infty}$ \[FMR].

For $J^{(2)}$ we have in turn
$$
\veqnalign{
J_{ij}^{(2)} =
&
2\suma2{k}{i-1}\suma2{l}{k}~
\qcomb[l-k-1/l]q^{(l-1)(k+1)}u_{j+k-l}~\d_q^{l}T^{-k}u_{i-k-1}
\cr
&
-~2\suma2{k}{i-1}\suma2{l}{j+k}\suma2{m}{i-k-1}~\qcomb[j-n-1/l]
\qcomb[n-m/i-k-m-1]
\cr
&~~~~~~~~~~~~~~~~~~~~~
q^{(l-1)(l-j+n+1)+(i+l-k-m-2)(i-k-n-1)}~u_m\d_q^{i+l-k-m-1}
T^{j+k-i-l+1}u_{j+k-l}
\cr
&
+~2\suma2{k}{i}\suma2{l}{k}~\qcomb[l-k-1/l]q^{(l-1)(k+1)}(q-1)~x~u_{j+k-l}
\d_q^{l}\,T^{-k}u_{i-k}
\cr
&
-~2\suma2{k}{i}\suma2{l}{j+k}\suma2{m}{i-k}~\qcomb[j-n-1/l]\qcomb[n-m/i-k-m]
\cr
&~~~~~~~~~~~~~~~~~~~~~
q^{(l-1)(l-j+n+1)+(i+l-k-m-1)(i-k-n)}(q-1)~x~
u_m\d_q^{i+l-k-m}T^{j+k-i-l}u_{j+k-l}
\cr
&
-~(1-q^{-1})~x~ u_i
u_j+\suma2{k}{i}\suma2{l}{j}~\qcomb[j-n-1/l]\qcomb[n-k/i-k]~
\cr
&~~~~~~~~~~~~~~~~~~~~~
q^{(l-1)(l-j+n+1)+(i+l-k-1)(i-n)}(q-1)~x~u_k
\d_q^{i+l-k}T^{j-i-l}u_{j-l}
\cr
&
+~\suma2{k}{i-1}\suma2{l}{j}~\qcomb[j-n-1/l]\qcomb[n-k/i-k-1]
\cr
&~~~~~~~~~~~~~~~~~~~~~q^{(l-1)(l-j+n+1)+(i+l-k-1)(i-n-1)}
\left(q^{n-i+1}-1\right)
{}~u_k \d_q^{i+l-k-1}T^{j-i-l+1}u_{j-l}
\cr
&+~\suma2{k}{i}~\qcomb[n-k/i-k]q^{(i-k)(i-n)}(1-q^{-1})~u_k
\d_q^{i-k}T^{n-i}~x~u_j
\cr
&
-~\suma2{k}{j}~\qcomb[j-n-1/k]q^{(k-1)(k-j+n+1)}(q-1)~x~u_i
\d_q^kT^{j-k-n}u_{j-k}
\(lagorda)
\cr}
$$
This expression reduces in the limit $q\to 1$ to the one of the
$W_{KP}$ algebra.
Contrarily to $J^{(1)}$, $J^{(2)}(X)$ does not stabilize the field $u_0$; \ie,
from \(elsegun) we see that the highest order of $J^{(2)}(X)$ is the same as that of
$L$. Therefore the constraint $u_0=1$ is second class. The same discussion that was
developed in the $T$ basis holds here {\sl mutatis mutandi}.  We will refrain from
giving the explicit form of the reduced Poisson brackets, whose computation follows
again the standard Dirac's recipe.

\subsection{Reductions. Where is $q$-$W_n$?}

Let us consider here the very important reductions of $q$-KP to
$q$-KdV. The expressions in \(elprimer) and \(elsegun) are perfectly consistent
when applied to purely $q$-differential operators $L= u_0\d_q^n+ u_1\d_q^{n-1}+...+
u_n$. The related algebras are simply obtained by restricting the
subindices of the fields appearing in \(mayorquene) \(jigene) and \(lagorda) to take
values in the range $i,j\in[0,n]$, and neglecting all other fields. In strict sense,
these algebras should be considered as deformations of $GD_n$, the
second Gelfand-Dickey bracket over the phase-space of Lax operators of the form
$L=\d^n+u_1\d^{n-1}+...+u_n$. Hence we will name them $q$-$GD_n$ algebras. 

An important point arises here: 
as compared with $GD_n$, $q$-$GD_n$ contain an
additional generator $u_0$. In the limit $q\to 1$ this field decouples
because $\lim_{q\to 1} J_{0j}= 0,~ \forall j$ and we may set $u_0=1$.
One could argue that in order to construct a true $q$-deformation
of $GD_n$ which involves exactly $n$ generators we should first reduce $u_0=0$ via
Dirac brackets. However the projection involved in the reduction {\sl is not a
continuous} step  and nothing guarantees that the resulting algebra will still
recover the desired limit when $q\to 1$.

Let us give an example of this phenomenon by considering the simplest Lax operator
$L=u_0\d_q+u_1$. The Poisson brackets for $u_0$ and $u_1$ generate $q$-$GD_1$,
whose brakets are given by
$$
\veqnalign{
J^{(2)}_{00} &= {q-1\ov 2q} u_0(T-T\inv)z \,u_0 \cr
J^{(2)}_{01} &= {1\ov 2q} u_0(T-T\inv)u_0 \cr
J^{(2)}_{10} &= {1\ov 2q} u_0(qT-q\inv T\inv)u_0 \cr
J^{(2)}_{11} &= {1\ov 2zq(q-1)} u_0(T-T\inv )u_0 
\(ceroyuno)\cr}
$$
which in the limit $q\to 1$ reproduce the free boson algebra $GD_1$ after $u_0$ is
set to 1,
\ie\  $J^{(2)}_{11}\to \d$, and $J^{(2)}_{0i}\to 0$. However if we insisted in
reducing
$u_0=1$ before taking the limit, the Dirac formula gives us a vanishing answer for
$\wt J^{(2)}_{11}$:
$$\veqnalign{
\wt J^{(2)}_{11} &= J^{(2)}_{11} - J^{(2)}_{10}(J^{(2)}_{00})\inv J^{(2)}_{01} \cr
&= u_0{1\ov 2q(q-1)}\left({1\ov z}(T-T\inv) - (qT-q\inv T\inv){1\ov z}\right)u_0
\cr  &=0  \(dacero)\cr}
$$    
One cannot cure this result by multiplying the starting brackets by global factors
of $(q-1)$, because the Dirac bracket is homogeneous under such rescalings. 
This vanishing result is also independent of any $q$-dependent redefinition
of the field $u_0$.

We expect that a similar discussion applies to the classical $W_n$ algebras although
we do not have a general proof. These algebras arise as hamiltonian reductions of
$GD_n$ where the generator $u_1$ is set to 0. The first generator, $u_2$,
closes a linear subalgebra which is non other than the ubiquitous Virasoro algebra.
It is in this sense that $W_n$ algebras are sometimes defined as (non-linear)
extensions of the Virasoro algebra.
 A continuous $q$-deformation of $W_n$ in $n-1$ fields $u_2,...,u_n$ would present
the same problems that we have exposed above in the case of $GD_n$. The naive
procedure, of starting from $q$-$GD_n$ and reducing $u_0=1$ and
$u_1=0$ may spoil the continuous correspondence with $W_n$ in the limit $q\to 1$.
We feel this is an important point that deserves further attention.

\section{Analitic continuation}

Notice that the expresion for the algebra $\qWKP^{(n)}$ as given in
\(lagorda) admits analitic continuation to complex values of $n=\alpha\in\comps$.
This happens in contrast with the first structure, given in \(mayorquene)\(jigene), 
where $n$ appears explicitely in the limits of sumatories. The best way to
understand this is by implementing the analytic continuation right from the
beginning. Actually the whole formalism is susceptible of such a continuation along
the lines advocated in ref. \[FMR] and \[KeZa]. Hence
$\qWKP^{(\alpha)}$ describes a {\it two- parameter} family of nonlinear
$\W_\infty$ type algebras.

There is an important technical question concerning the triviality of
such deformation parameters, \ie, whether the algebras for differente pairs
$(q,\alpha)$ and $(q',\alpha')$ are isomorphic or not. At least in the
continuum case $q=1$ we know positively that $\alpha$ represents a non-trivial
deformation parameter \[Fernando].

The second issue we intend to address in this concern is the possibility of
connecting the linear and quadratic structures by a suitable contraction of
the parameter $\alpha$. In \[FMR] the limit $\alpha \to 0 $ was shown to
yield an extension of the linear algebra $\W_{1+\infty}$ by means of the
Khesin-Kravchenko cocycle \[KK]. This fact was also understood in \[BaKeKi] and
in \[KeZa] from a Poisson-Lie group theoretical point of view.

In more concrete terms, let us introduce
a second parameter $\beta$ and define $L^{(\alpha,\beta)}\in
\M^{(\alpha)}_{\d_q}$ such
that  
$$
L^{(\alpha,\beta)}= \beta \d^{\alpha}_q+\sum_{j=0}^\infty u_j(z)\d^{\alpha - j}_q
\equiv \beta \d^{\alpha}_q + L^{(\alpha)}.
\(separacion)
$$
Correspondingly, the 1-forms $X$ will look as
$$
X= \sum_{j=0}^\infty \d^{j-\alpha-1}_q ~ T~ x_j.
$$
We will be
interested in the following ``scaling" limit  in which $\alpha$ tends to 0 and
$\beta$ to $\infty$ in such a way that $\alpha\beta = c$ a finite constant.
It will be convenient to normalize $J^{(2)}$ in the following form
$$\veqnalign{
J^{(2)}_{L^{(\alpha,\beta)}}(z) =& {1\over \beta} \left\{
L^{(\alpha,\beta)}(X L^{(\alpha,\beta)})_{\geq 0} - (L^{(\alpha,\beta)}X)_{\geq
0}
L^{(\alpha,\beta)}\right.
\cr   & \left. +\med L^{(\alpha,\beta)}~\Omega (\comm{L^{(\alpha,\beta)}}{X}) +
\med
\Omega(\comm{L^{(\alpha,\beta)}}{X})~L^{(\alpha,\beta)}
\right\} \cr}
$$
Plugging \(separacion) in this expression we may first gather all the terms
quadratic in $\beta \d_q^{\alpha}$:
$$
\beta \left( \d_q^{\alpha} (X \d_q^{\alpha})_{\geq 0} -
(\d_q^{\alpha} X )_{\geq 0}  \d_q^{\alpha}
+ \med\d_q^{\alpha} (\Omega\comm{\d_q^{\alpha}}{X})
+ \med \Omega(\comm{\d_q^{\alpha}}{X})  \d_q^{\alpha} \right)
$$
Expanding $ \d_q^{\alpha} = 1+\alpha\log\d_q+ O(\alpha^2)$
the surviving terms in the desired limit yield
$$
c\comm{\log\d_q}{X_{\geq 0}} -c \comm{\log\d_q}{X}_{\geq 0}
+c\Omega(\comm{\log\d_q}{X})
= c~\comm{\log\d_q}{X_{\geq 0}}_{\leq -1}  +c~\Omega(\comm{\log\d_q}{X})
$$
In the linear terms the $\beta$ dependence cancels out and
we obtain
$$
\comm{ L^{(0)}}{X_{\geq 0}} - \comm{L^{(0)}}{X}_{\geq 0}
+\Omega(\comm{L^{(\alpha)}}{X})
= \comm{L^{(0)}}{X_{\geq 0}}_{\leq -1}
+\Omega(\comm{L^{(0)}}{X})
$$
In summary, the limiting hamiltonian structure yields
$$
J_{1+\infty;q}(X) = \comm{c\log\d_q + L^{(0)}}{X_{\geq 0}}_{\leq -1}
            +\Omega(\comm{c\log\d_q + L^{(0)}}{X}) \(qwinfty)
$$
Consistency of $J(X)$ as a tangent map demands that $L$  be of the form
$$
L= \log \d_q + u_0 + \sum_{i=1}^\infty u_i \d^{-i}
$$
The expression
 $\log \d_q$ has to be understood as an outer automorphism of the
Lie algebra of
$q$-pseudodifferential symbols. Its action can be defined and computed as a
limit:
$$\veqnalign{
\comm{\log\d_q}{f\d_q^p} &= \lim_{\alpha\to 0}~{1\ov \alpha}(\d_q^\alpha
{}~f\d_q^p~ \d_q^{-\alpha} -f\d_q^p)
\cr
&
 = \log q f'\,\d_q^p -\sum_{k\geq 1}
{\log q\ov (q-1)} \qcomb[-1/k] q^k~\tau^{-k}(\d_q^k
f)
{}~\d_q^{p-k}  \cr}
$$
The notation in \(qwinfty) intends to make explicit that this algebra is a
$q$-deformation of the centrally extended $W_{1+\infty}$ algebra. We will not
write
down explicitely the Poisson brackets here. They agree with the ones given
in  \(mayorquene) and  \(jigene) with $n=0$ except for the central terms, which
are the only ones that acquire corrections  proportional to $\log q$.

\section{Conclusions and outlook}
The picture of an atlas of $\W$ algebras is slowly emerging.
 In this landscape, $W_\infty$ algebras provide natural landmarks and, among
them, the algebra $\WKP$ is a cornerstone. In ref.\[atlas] this algebra was
shown to be related with a large amount of the known classical W-type algebras
by continuous deformation or truncation.  The main result of the present paper
is that a lot of points in that atlas admit yet another deformation,
parameterized by $q$. Of particular importance are $q$-$\WKP^{(\alpha)}$,
$q$-$GD_n$ and $q$-$W_{1+\infty}$.

It has been amusing to observe how many structures that worked fine for the algebra
of pseudodifferential operators, are robust enough to resist their implementation
in the algebra of $q$-deformed pseudodifferential operators, as well. It certainly
points out that perhaps other well known results could be exported. To be more
precise, we think about issues like the dressing transformation,
the embedding of the Lie algebra of differential operators into $W_n$ \[RadulFigRa]
or  the Kupershmidt-Wilson-Yu theorem \[KuWiYu]. In fact, concerning this last
important theorem, a straightforward implementation of the proof given in
\[KPMiura] for a quadratic structure of the form \(quadratic) 
works fine in the case of an isotropic splitting. 
 This requirement is only fulfilled in the present work for the splitting 
$\sigma=0$ and hence there is a q-deformed version of the Kupershmidt-Wilson-Yu
theorem in this case. For $\sigma=\pm1$ we have not been able to establish a
similar result. In this respect we should mention that a proposal for a q-deformed
Miura transformation has appeared in \[FResh]. Its connection to some peculiar 
way to factorize the Lax operator has been addressed in \[Frenkel].

We should emphasize the existence of three consistent
splittings ($\sigma=0,\pm 1$), for the algebra $\pdo_q$. They all yield
integrable hamiltonian systems and thereafter $\W$ type algebras.
In references \[Kup]\[OeStra] the m-KdV hierarchy was investigated in the
scalar Lax formalism. It was recognized that this system is related to a nonstandard
splitting of the algebra of ordinary \pdo. Indeed $L = L_{\geq k} + L_{< k}$ yields
consistent subalgebras for $k=0,1$ and
$2$. It would be interesting to find out whether a possible $q$-deformation of these
non-standard splittings could be related to the cases $\sigma=0$ and $\sigma = +1$
in this paper.

The connection of the KP hierarchy with the Toda lattice hierarchy is a subject
of recent
interest which has received the attention of different groups
\[Gieseker]\[Bonora]\[Aratyn]. We believe that our approach is substantially
different
to these and closer, at least in spirit, to the lattice deformation of
\[Babelon]. We expect that the powerfull techniques that have been used in this
paper
can be implemented also in the context of the Calogero-Sutherland model,
especially in the
formulation that makes use of the exchange operators
\[Calogero].

\ack
We are indebted to Boris Khesin for his useful comments and for bringing
several references in \[FResh]\[QuantW] to our attention, and to E. Frenkel
for clarification about his work. We also want to thank warmly Eduardo Ramos
for his reading and  comments on the manuscript. This work has been partially
supported by the DGICYT under contract PB/93-0344.

\refsout
\bye